%% file: main.tex
\documentclass[11pt,a4paper]{amsart}
\input{preamble.tex}

\title[Neural Network Acceleration]{Neural Network Acceleration of Iterative Methods for nonlinear Schrödinger eigenvalue problems}

\author[D.~Peterseim]{Daniel Peterseim}
\author[J.~F.~Pietschmann]{Jan-F.~Pietschmann}
\author[J.~Püschel]{Jonas Püschel}
\author[K.~Rue\ss]{Kilian Rue\ss}
\address[D.~Peterseim, J.-F.~Pietschmann]{Institute of Mathematics \& Centre for 
Advanced Analytics and Predictive Sciences (CAAPS), University of 
Augsburg, Universit\"atsstra{\ss}e~12a, 86159 Augsburg, Germany}
\email{daniel.peterseim@uni-a.de}
\email{jan-f.pietschmann@uni-a.de}
\address[J.~Püschel, K.~Rue\ss]{Institute of Mathematics, University of Augsburg, 
Universit\"atsstra{\ss}e~12a, 86159 Augsburg, Germany}
\email{jonas.pueschel@uni-a.de}
\email{kilian.ruess@uni-a.de}

\thanks{The work of D.~Peterseim is part of a project that has received funding from the European Research Council (ERC) under the European Union's Horizon 2020 research and innovation programme (Grant agreement No.~865751 -- RandomMultiScales).}

\begin{document}

\begin{abstract}
We present a novel approach to accelerate iterative methods to solve nonlinear Schrödinger eigenvalue problems using neural networks. Nonlinear eigenvector problems are fundamental in quantum mechanics and other fields, yet conventional solvers often suffer from slow convergence in extreme parameter regimes, as exemplified by the rotating Bose-Einstein condensate (BEC) problem. Our method uses a neural network to predict and refine solution trajectories, leveraging knowledge from previous simulations to improve convergence speed and accuracy. Numerical experiments demonstrate significant speed-up over classical solvers, highlighting both the strengths and limitations of the approach. 
\end{abstract}

\maketitle
{\tiny {\bf Key words:} nonlinear Schrödinger, Gross-Pitaevskii, energy-adaptive Riemmanian conjugate gradient descent, neural network, U-net}

{\tiny {\bf AMS subject classifications.} 
65H17, 65N25, 81Q10, 35Q55, 58C40
}


\section{Introduction}

Nonlinear eigenvalue problems (NLEVPs) with eigenvector nonlinearity, that is, where the eigenvector appears nonlinearly in the operator, frequently appear in computational physics and chemistry. Classical examples in Schrödinger‐type settings include the calculation and prediction of the properties of molecules and solid-state materials via Hartree-Fock \cite{Slater1951} and Kohn-Sham \cite{KohS65} equations of electronic-structure theory. Here, the unknown electron density enters the potential nonlinearly. Another equally important case is the Gross-Pitaevskii equation (GPE), which describes Bose-Einstein condensates (BECs) \cite{PhysRevLett.78.586,SteISMCK98} in the mean-field approximation \cite{LieSY00}.

Classical numerical methods for ground‑state computations in nonlinear eigenvalue problems (NLEVPs) have reached a high level of sophistication. They either operate directly on the NLEVP, such as self‑consistent field (SCF) iterations \cite{DioC07,CanKL21}, or exploit that, in the applications above, the NLEVP is the Euler–Lagrange equation of an energy minimization problem subject to normalization (or orthonormality) constraints. Established methods of the latter class include discrete normalized gradient flows \cite{BaoD04,BaoCL06}, (projected) Sobolev gradient methods \cite{GarP01,KaE10,DanK10,HenP20,ChenLLZ24}, the $J$‑method \cite{JarKM14,AltHP21}, Riemannian optimization methods in discrete and continuous settings \cite{AntLT17,DanP17,AltPS22,PPS25}, and Newton‑type approaches \cite{BaoT03,DuL22,WuWB17,AltPS24}. For the Gross–Pitaevskii equation (GPE) in particular, an extensive overview is provided in \cite{HenJ23}.

Yet classical solvers for these NLEVPs can become prohibitively expensive in challenging parameter regimes that often underlie the most interesting physics. While electronic‑structure calculations primarily struggle with accurately discretising singular, long‑range Coulomb interactions, in the Gross–Pitaevskii equation (GPE) the characteristic cubic nonlinearity itself dominates and frequently hampers conventional self‑consistent field (SCF) iterations. This difficulty is further amplified in rotating condensates with high angular momentum, where dense vortex lattices make the energy landscape highly non‑convex. Even state-of-the-art Riemannian optimization schemes may require several thousand or even ten thousand iterations and can stagnate in non-global stationary states (see also Table~\ref{tab:convergence_earcg_1}  below). For these reasons, we make use of the rotating GPE in extreme parameter regimes as a demanding yet representative benchmark for NLEVP solvers. 

Recent data‐driven advances have demonstrated significant potential in enhancing the performance of numerical algorithms. For example, machine learning techniques have been used to speed up iterative algorithms for convex optimization \cite{DBLP:journals/corr/abs-2107-10254} (coined neural fixed-point acceleration) and in the context of inverse problems \cite{adler2017solving}. However, their extension to NLEVPs involving nonconvex manifold constraints has remained unexplored.

Here, we present a neural‐network accelerator for Riemannian optimization of nonlinear Schrödinger eigenvalue problems. By embedding a lightweight U-Net based neural network into mid‐iteration steps of the energy‐adaptive Riemannian conjugate‐gradient (EARCG) solver \cite{PPS25,HenP20,AltPS22}, our approach circumvents stagnation and on average reduces iteration counts by $22\%$ and wall‐clock time by $14.5\%$ in rotating GPE benchmarks. The neural network is trained offline on representative discrete solution trajectories of the EARCG.

Unlike end-to-end neural solvers that replace classical algorithms with black-box models \cite{BAO2025113486,PhysRevResearch.7.013332}, our hybrid approach uses the network strictly as an accelerator, retaining the solver’s favorable robustness and convergence properties. It also differs from methods that use Riemannian optimization to train neural networks. Because many other NLEVPs, such as multicomponent BEC systems \cite{altmann2025riemannianoptimisationmethodsground} and Kohn–Sham–type electronic-structure models \cite{AltPS22,AltPS24,PPS25}, employ similar Riemannian solvers, our workflow can be transferred to those settings with minimal adaptation, promising an order-of-magnitude speedup across a broad class of applications.

The remainder of the paper is organized as follows. Section~\ref{sec:model} introduces the rotating GPE and its physical parameters. Section~\ref{sec:riemann} revisits the EARCG algorithm within a larger context of Riemannian optimization methods. Section~\ref{sec:nna} details the neural-network accelerator, including data generation, architecture, and training. Section~\ref{sec:numexp} reports extensive numerical experiments in critical parameter regimes where vortex lattices emerge, and Section~\ref{sec:conclusion} concludes with perspectives.

\section{Model Benchmark Problem}\label{sec:model}

The rotating Gross-Pitaevskii equation will serve as a model problem for our approach, and we briefly reiterate its formulation together with a list of relevant parameters below. 
For a more in-depth discussion of the Gross-Pitaevskii models, the interested reader is referred to \cite{BaoCai2013, HenJ23}.

Let $\mathcal{D} = [-\tfrac{a}{2},\tfrac{a}{2}]^2 \subset\IR^2$ be a sufficiently large square domain of width $a > 0$. Here and throughout, we choose a square for simplicity, but any convex bounded domain in one, two, or three dimensions would be equally feasible. On $\mathcal D$, we consider the Hilbert space $\Lsp \coloneqq L^2(\mathcal{D},\C)$ and the Sobolev space $\Hsp \coloneqq H_0^1(\mathcal{D}, \C)$ as well as its dual space $\Hspi \coloneqq H^{-1}(\mathcal{D}, \C)$. The former is endowed with the real inner product
\[
(v, w)_{\Lsp} 
:= \re \int_{\mathcal{D}} \overline{v} \, w \,\dx  
\]  
which also corresponds to the dual-evaluation on $\Hspi \times \Hsp$ as a consequence of the Gelfand triple structure $\Hsp \subset \Lsp \subset \Hspi$.

Our main object of study is the Gross-Pitaevskii (GP) energy functional $\calE\colon \Hsp \rightarrow \IR$ defined by 
\begin{equation}\label{eq:energy_short}
\calE(\varphi) := \int_\mathcal{D} \frac{1}{2} |\nabla \varphi|^2
+ \frac{1}{2}V|\varphi|^2  + \frac{1}{2} \omega ( \overline{\varphi}  L_z  \varphi) + \frac{1}{4}\kappa |\varphi|^4 \,\mathrm{d}x
\end{equation}
for any $\varphi\in \Hsp$, where
\begin{align*}
    V(x) = v_1 x_1^2 + v_2 x_2^2,\quad v_1,v_2>0
\end{align*}
is the harmonic trapping potential,
\[
L_z \;=\;-\,i\bigl(x_1\partial_{x_2}-x_2\partial_{x_1}\bigr)
\]
the \(z\)-component of angular momentum (with rotation speed \(\omega>0\)),
and \(\kappa>0\) the interaction strength.  Since \(V(x)\to+\infty\) as \(|x|\to\infty\),
ground states decay exponentially, so imposing zero Dirichlet boundary conditions on \(\partial\mathcal D\)
introduces only negligible truncation error provided \(a\) is chosen sufficiently large relative to $v_1$, $v_2$, and $\kappa$.

Our aim is to compute a global minimizer of $\calE$, a~so-called {\em ground state}, on the set 
\begin{equation}\label{eq:sphere}
	\SH = \left\{\varphi \in \Hsp \mid \|\varphi\|_{\Lsp} = 1\right\},
\end{equation}
which carries the structure of a Riemannian manifold. The normalization in $\Lsp$ ensures that the densities $|\varphi|^2$ are probability densities.
This yields the constrained energy minimization problem 
\begin{equation}\label{eq:min}
	\min_{\varphi\in\SH} \calE(\varphi)
\end{equation} 
on the Riemannian manifold $\SH$.
Define the bilinear forms
\begin{equation}\label{eq:a}
    \begin{aligned}
         a_\varphi(v, w) = \re \Big(\int_\mathcal{D} \nabla \overline{v}\cdot \nabla w
+ V\overline{v}w  + \omega ( \overline{v}  L_z  w) + |\varphi|^2 \kappa ( \overline{v}w) \,\mathrm{d}x\Big)
    \end{aligned}
\end{equation}
on $\Hsp\times \Hsp$, which depend on the density $|\varphi|^2=\overline{\varphi} \varphi$ of the given state $\varphi \in \Hsp$. The directional derivative $\Drm \calE(\varphi)[v]$ of $\,\calE$ at $\varphi \in \SH$ is then compactly written as  
\begin{equation*}\label{eq:D_E}
    \Drm \calE(\varphi)[v] = a_\varphi(\varphi, v).
\end{equation*}
Given $a_\varphi(\cdot,\cdot)$, we also define operators $\A_\varphi: \Hsp \to \Hspi$ by
\begin{equation}\label{eq:A}
    \langle \A_\varphi(v),w \rangle_{\Lsp} = a_\varphi(v,w) 
\end{equation}
for all $v,w \in \Hsp$. By definition, these operators are isomorphisms and they are Hermitian with respect to $\langle \cdot, \cdot \rangle_\Lsp$. A strong form representation reads
\begin{equation}\label{eq:A_explicit}
    \A_\varphi = -\Delta + V + \omega \mathcal L_3 + |\varphi|^2 \kappa.
\end{equation}
In the context of non-linear eigenvalue problems, $\A$ is often referred to as the \emph{Hamiltonian operator}.
By the method of Lagrange multipliers, a necessary and sufficient condition for $\varphi$ to be a critical point of $\E$ on $\SH$ is the existence of a real-valued Lagrange multiplier $\lambda \in \mathbb R$ such that 
\begin{equation}\label{eq:evp}
    \A_\varphi(\varphi) = \lambda \varphi.
\end{equation}
Thus, ground states of $\E$ can either be calculated using energy minimization of $\E$ on $\SH$ or via solving the self-consistent eigenvalue problem \eqref{eq:evp}. Heuristically, small eigenvalues correspond to small energies. However, the equivalence of the smallest eigenvalue with the lowest energy is proven only in special cases, such as when there is no rotation ($\omega = 0$), as discussed in \cite{CanCM10, HLP24}. However, this equivalence does not hold in general cases involving rotation, see \cite{AltHP21}.

\begin{figure}[h]
    \centering
    \includegraphics[width=0.95\linewidth]{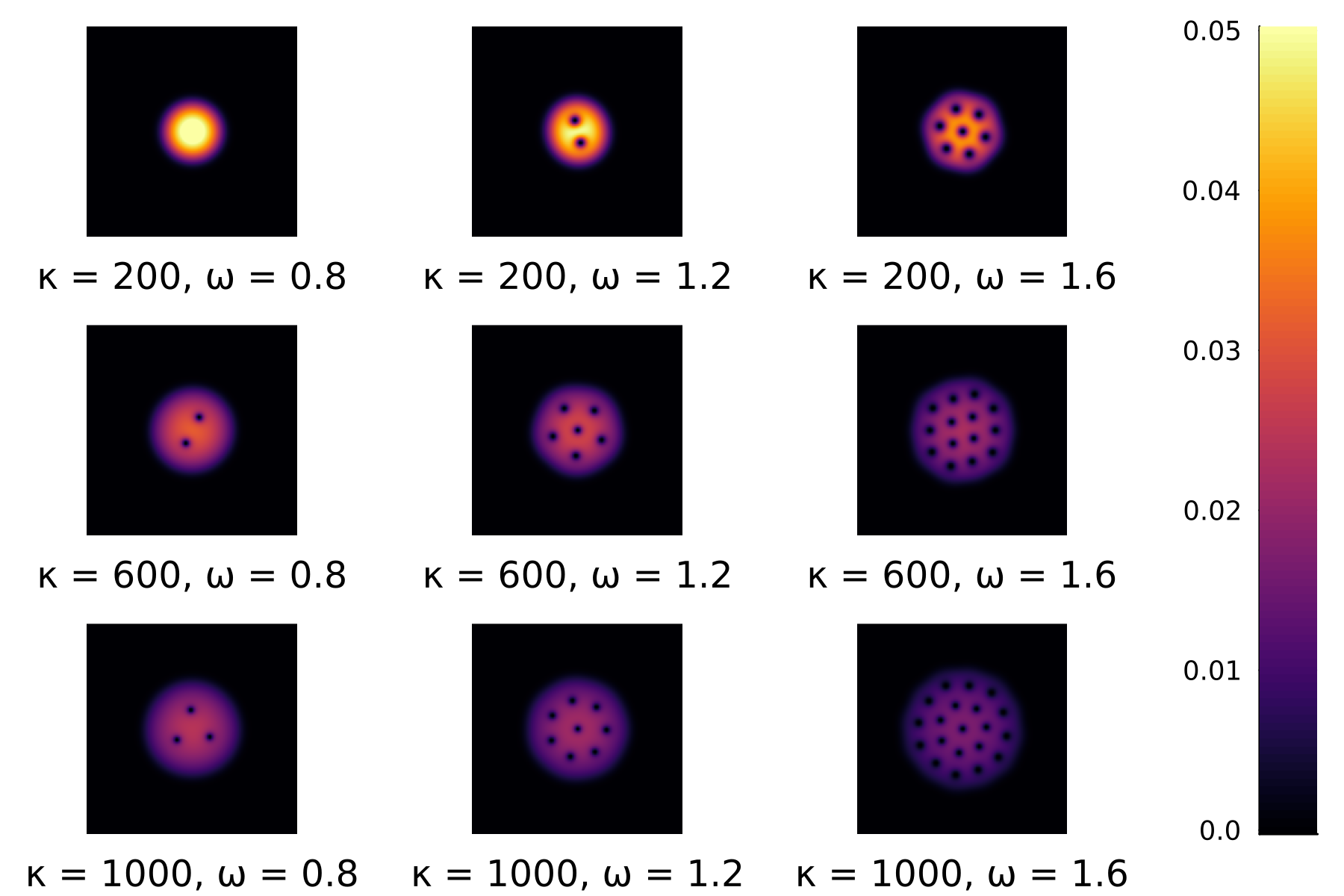}
    \caption{Local minimizer densities for $a =20,\;v = (1,1)^T$}
    \label{fig:ground_states_1}
\end{figure}

\begin{figure}[h]
    \centering
    \includegraphics[width=0.95\linewidth]{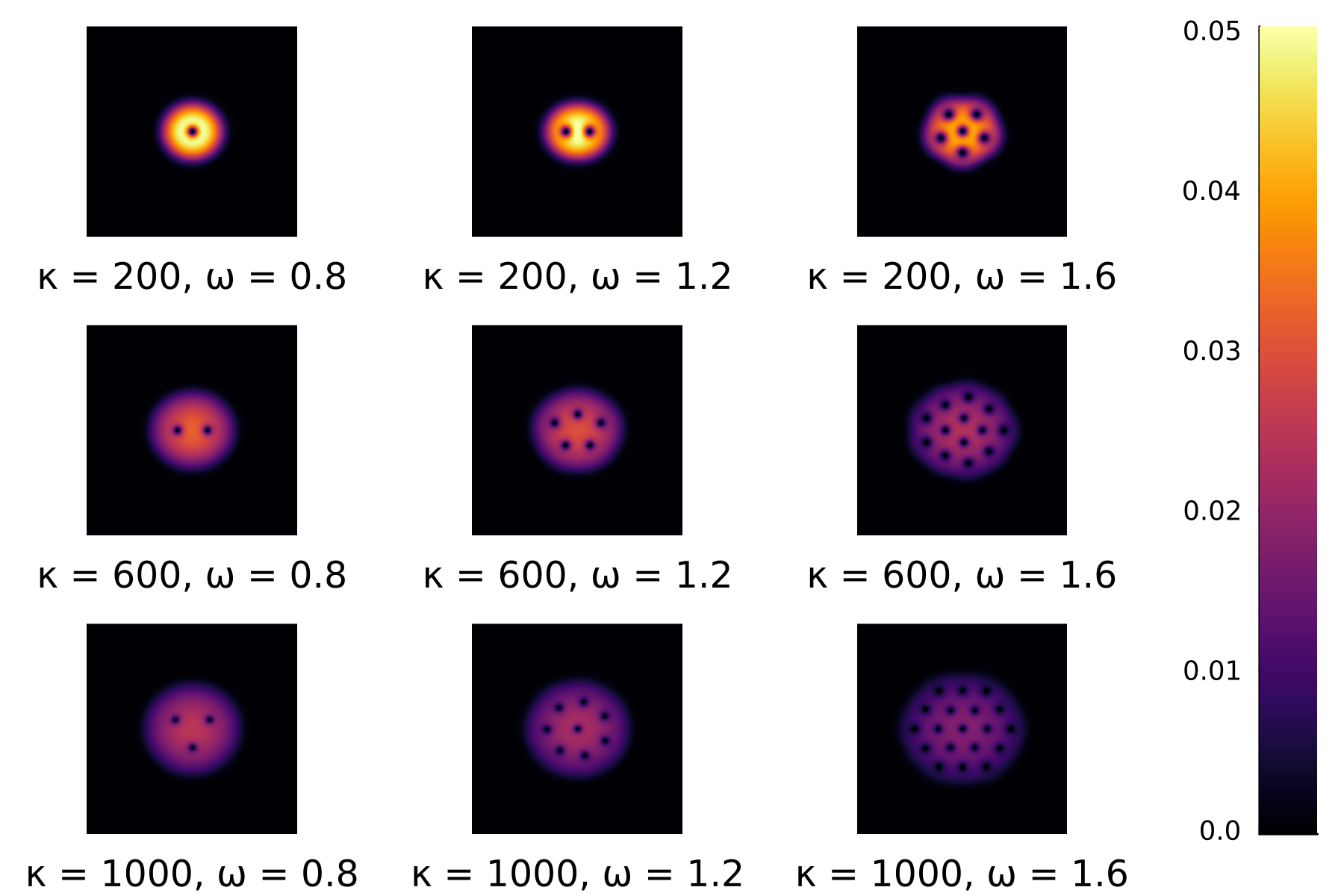}
    \caption{Local minimizer densities for $a =20,\;v = (1.1,1)^T$}
    \label{fig:ground_states_2}
\end{figure}

\begin{figure}[h]
    \centering
    \includegraphics[width=0.95\linewidth]{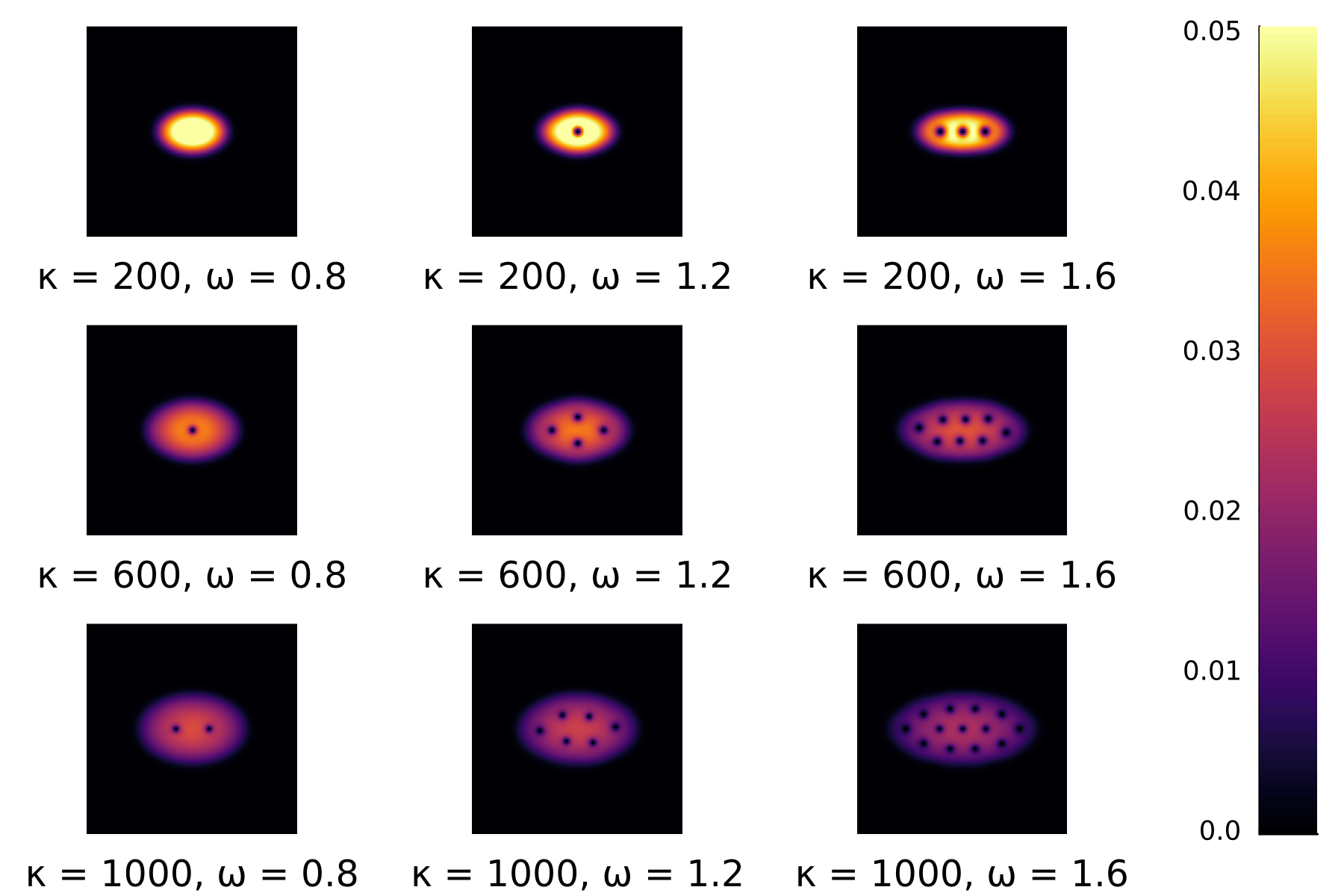}
    \caption{Local minimizer densities for $a =20,\;v = (2,1)^T$}
    \label{fig:ground_states_3}
\end{figure}

We briefly discuss the effects that the parameters introduced in the GP model above have on the ground state.
\begin{itemize}
\item $a$: The size of the cell does in general not affect the state of its density, as long as the cell is sufficiently large.
\item $v_i$: The amplitude of the confining potential in $x_i$-direction. Larger $v_i$ corresponds to more concentrated density in direction $x_i$. For $v_1 \neq v_2$, the confining is non-symmetric, i.e.\ the ground state is generally not invariant under rotation.
\item $\omega$: The amplitude of the rotating magnetic field. Stronger magnetic fields generate more vortices.
\item $\kappa$: The multiplicative coefficient of the quartic variational term (representing the cubic nonlinearity), which penalizes high densities. Thus, larger values of $\kappa$ lead to a larger spread of the density in the ground state, since the repulsion is stronger and thus wider mass spread reduces the potential energy. This also promotes vortex formation.
\end{itemize}
Figures \ref{fig:ground_states_1}--\ref{fig:ground_states_3} depict local minimizers of $\calE$ that illustrate the previously described effects of the different parameters.

\section{Energy-adaptive Riemannian conjugate gradient descent}\label{sec:riemann}

Riemannian optimization is one well-established methodology for solving \eqref{eq:min}.  Here, as a prime representative of more general Riemannian methods, we introduce the Energy-Adaptive Riemannian Conjugate Gradient (EARCG) method; later we will also highlight its limitations in extreme parameter regimes. For a more detailed description of the method on the more general Stiefel manifold as well as some numerical improvements, the interested reader is referred to \cite{PPS25}.

Recall the coercive, bounded, symmetric bilinear form \(a_\varphi\) from \eqref{eq:a}, for $\varphi \in \SH$, and the respective operator $\A_\varphi$ from \eqref{eq:A_explicit}. 
Then the \emph{energy-adaptive metric} at $\varphi \in \SH$ is given as
\begin{equation}\label{eq:ea_metric}
    g^a_\varphi(v,w) = a_\varphi(v,w),\quad v,w \in \Tan_\varphi \SH,
\end{equation}
where $\Tan_\varphi \SH = \{ v \in \Hsp \mid \langle\varphi, v \rangle_{\Lsp} = 0\}$ is the tangent space of $\SH$ at $\varphi$. This choice of inner product results in the  \emph{energy‐adaptive gradient} of \(\E\) at \(\varphi\) being
\begin{equation}\label{eq:ea_grad_E}
\grad_{\!a,\,\varphi}\E(\varphi)
=\varphi \;-\;\frac{\A_\varphi^{-1}\varphi}
{\langle \A_\varphi^{-1}\varphi,\;\varphi\rangle_{\Lsp}}.
\end{equation}
We combine this with the normalization retraction
\(\R^{\rm norm}_\varphi(v)=(\varphi+v)/\|\varphi+v\|_{\Lsp}\)
and its differentiated retraction vector transport
\(\T^{\rm norm}_v(w)
=D(\,u\mapsto u/\|u\|_{\Lsp}\,)\bigl[\varphi+v\bigr][\,w\,]\) and use the convergence criterion $\|g^{(k)}\|_{g^a_\varphi} < tol$ where
\begin{equation}\label{eq:ea_norm}
    \|\cdot\|_{g^a_\varphi} = \sqrt{g^{a}_\varphi (\cdot, \cdot)}
\end{equation}
 is the norm on $\Tan_\varphi\SH$ induced by the metric $g^a_\varphi$, which we will refer to as the \emph{energy norm}. Using the conjugate gradient (CG) parameter and step size strategy from \cite{PPS25}, we construct a Riemannian CG scheme, which we will refer to as the energy-adaptive Riemannian conjugate gradient (EARCG) method. It is summarized in Algorithm~\ref{alg:earcg}.

\begin{algorithm2e}[th]
\SetKwInOut{Input}{input}
\SetAlgoLined
\Input{initial guess $\varphi^{(0)} \in \SH$, tolerance $tol$}
\BlankLine
$g^{(0)}= \varphi^{(0)}- \A_{\varphi^{(0)}}^{-1}(\varphi^{(0)})\langle \varphi^{(0)}, \A_{\varphi^{(0)}}^{-1}(\varphi^{(0)}) \rangle_{\Lsp}^{-1}$\;
$\eta^{(0)} = - g^{(0)}$ \;
$ \tau^{(-1)} = 1$ \;
 \For{$k=0, 1, 2, \dots$ until $\|g^{(k)}\|_{g^a_{\varphi^{(k)}}} < tol$}{
    compute the step size $\tau^{(k)}$ with \cite[Algorithm 2]{PPS25} and initial guess $\tau^{(k-1)}$ \;
    $\varphi^{(k+1)} = \R^{\rm norm}_{\varphi^{(k)}}(\tau^{(k)} \eta^{(k)})$\; 
    $g^{(k+1)}= \varphi^{(k+1)}- \A_{\varphi^{(k+1)}}^{-1}(\varphi^{(k+1)})\langle \varphi^{(k+1)}, \A_{\varphi^{(k+1)}}^{-1}(\varphi^{(k+1)}) \rangle_{\Lsp}^{-1}$\;
    $\beta^{(k+1)} = \max\bigg\{0, \min\Big\{\frac{a_{\varphi^{(k+1)}}( g^{(k+1)}, g^{(k+1)} )}{a_{\varphi^{(k)}}( g^{(k)}, g^{(k)} )}, \frac{a_{\varphi^{(k+1)}}( g^{(k+1)}, g^{(k+1)} - \T^{\rm norm}_{\tau^{(k)} \eta^{(k)}} (g^{(k)}) )}{a_{\varphi^{(k)}}( g^{(k)}, g^{(k)} )}\Big\}\bigg\}$\;
    $\eta^{(k+1)} = - g^{(k+1)} + \beta^{(k+1)} \T^{\rm norm}_{\tau^{(k)} \eta^{(k)}} (\eta^{(k)})$\; 
    
    $\tau^{(k+1)} = \tau^{(k)}$\;
 }
 \KwRet{$\varphi^{(k+1)}$}
 \caption{EARCG for Gross--Pitaevskii}
     \label{alg:earcg}
\end{algorithm2e}

The principal challenge when minimizing the Gross–Pitaevskii energy with a rotational term via iterative methods is their slow convergence. Even the EARCG method, which typically requires far fewer iterations than methods based on the $H^1$ or $L^2$ metrics, can still demand a very large number of steps to drive the energy‐norm of the gradient below a given tolerance.  Table~\ref{tab:convergence_earcg_1} reports the iteration counts required for EARCG to reach $\|g\|_{a,\varphi}<10^{-8}$ for the ground‐state solutions shown in Figures \ref{fig:ground_states_1}–\ref{fig:ground_states_3}.  (We emphasize that these are single runs; iteration counts can fluctuate with different randomized initial guesses or when converging to different local minima.)  

Two factors dominate the convergence rate: the number and arrangement of vortices (larger values of $\kappa$ and $\omega$ generally decelerate convergence) and any trap asymmetry, $v=(v_1,v_2)$ with $v_1\approx v_2$, which induces competing directional biases that neither fully align nor completely dominate the vortex lattice, as shown in Figure \ref{fig:ground_states_2}. 
To illustrate this, consider $a=20$, $v=(1.1,1)$, $\omega=1.4$, $\kappa=1000$.  

\input{figs/heatmaps_convergence_classical}

As shown in Figure \ref{fig:slow-convergence}, the EARCG iteration rapidly discovers the correct vortex pattern (by about the \(500\)th step) but then requires almost $12000$ further iterations to rotate that pattern into alignment with the trap and drive down the potential energy, finally converging at iteration $12140$.  By contrast, for the choice $v=(1,1)$ the method terminates after only $1135$ iterations (no rotation is needed since the ground state is rotationally invariant), and for the strongly anisotropic case $v=(2,1)$ it takes only $904$ iterations (the initially found configuration is already correctly oriented). This significant difference underscores how a slight trap asymmetry can dramatically slow convergence by forcing an expensive realignment phase.

\begin{table}[H]
    \centering
   \begin{tabular}{|c||c|c|c||c|c|c||c|c|c|}\hline
\diagbox{$\kappa$}{$\omega$}
    & \multicolumn{3}{c||}{\makebox[3em]{$0.8$}}&\multicolumn{3}{c||}{\makebox[3em]{$1.2$}}&\multicolumn{3}{c|}{\makebox[3em]{$1.6$}} \\\hline\hline
         $200$   & 441 & 306 & 111 & 204 & 359 & 292 & 334 & 1008 & 412 \\ \hline
         $600$   & 529 & 832 & 370 & 493 & \textbf{13164}  & 502 & 792 & 2816  & 673 \\ \hline
         $1000$  & 540 & \textbf{7075}  & 500 &  1096 & \textbf{15678}  & 917 & 1169 & \textbf{9614}  & 1299 \\ \hline
\end{tabular}
    \caption{Number of iterations of EARCG to reach \texttt{tol = 1e-8} for $a =20$ with $v = (1,1)^T$ (left), $v = (1.1,1)^T$ (middle) and $v = (2,1)^T$ (right)}
    \label{tab:convergence_earcg_1}
\end{table}

\section{Neural Network Acceleration}\label{sec:nna}
We now present a neural network enhanced version of the previously described EARCG method. First, we introduce the general acceleration strategy. We then discuss data generation, neural network architecture, and training as well as the specific acceleration strategy.

\subsection{Hybrid computational framework}

We fix the parameters \(0<\epsilon_2 < \epsilon_1^{\min}<\epsilon_1^{\max}\), $n_e \in \mathbb{N}$ as well as  \(e_0>0\). Then, our acceleration algorithm consists of three phases:
\begin{enumerate}\setlength{\itemsep}{2pt}
    \item 
\textbf{Initial EARCG phase}: We execute the EARCG algorithm until the energy norm \eqref{eq:ea_norm} of the gradient falls below the threshold $\epsilon_1^{\max}$. 
    \item 
\textbf{Neural Acceleration}: If the energy-adaptive norm of the gradient lies within the interval $[\epsilon_1^{\rm min}, \epsilon_1^{\rm max}]$, we apply the neural network every $n_e$-th iteration. Its output is accepted only if the normalization error, i.e.~the absolute difference of the squared $L^2$ norm from one, is smaller than $e_0$. If the energy-adaptive norm of the gradient reaches $\epsilon_1^{\rm min}$ without any neural network output having been accepted, neural network acceleration is enforced once. 
    \item
\textbf{Final EARCG Phase}: 
After neural network acceleration, we apply a final phase of EARCG until the gradient reaches the desired convergence threshold $\epsilon_2$
    in the energy norm.
\end{enumerate}  

We will explain the neural network acceleration in detail below.

\subsection{Training data}\label{sec:training_data}
The training dataset is constructed by running EARCG on random initial data $\varphi_0 \in \SH$ until the tolerance reaches the threshold $\epsilon_2$. For random generation of $\varphi_0$, first every component (i.e. coefficients of the plane wave basis used by DFTK, see Section~\ref{sec:numexp} below) is sampled from a normal distribution with mean $0$ and standard deviation $1$ and then the resulting vector is $\Lsp$-normalized.  
We then use intermediate iterates $\varphi_j$ and energy-adaptive Riemannian gradients $g_j$, chosen at $20$ log-equidistantly spaced tolerances $\tilde \epsilon_j \in [\epsilon_1^{\min}, \epsilon_1^{\max}]$, resulting in 20 pairs of training data per run. Each training data point then consists of a pair $((\varphi_j, g_j), \varphi^*)$, where $\varphi^*$ denotes the converged solution of the EARCG method.

The orbital $\varphi$ generated by DFTK requires preprocessing to match the input requirements of our image-based neural network. We apply an inverse Fourier transform to $\varphi$, resulting in a complex-valued matrix, which we subsequently represent as a $(n, n, 2)$ real-valued tensor, where the two channels correspond to the real and imaginary components, respectively.

\subsection{Neural network architecture and training}
Our implementation employs a modified U-Net architecture \cite{unet2015}, originally developed for image segmentation tasks, containing $31$ million trainable parameters. Padding is used in $3\times3$ convolutions to prevent loss of border pixels and thus merging outputs does not require cropping. The network accepts a $(n, n, 4)$-tensor as input, comprising the (real and imaginary part of the) intermediate solution state and its energy-adaptive gradient. It produces a $(n, n, 2)$ tensor representing (again the real and imaginary part of) the enhanced state. The network architecture is visualized in Figure \ref{tikz:network}.
The model was trained for $100$ epochs using the Adam optimizer with an initial learning rate of $10^{-4}$. 
We use the loss function $L(\hat \varphi, \varphi) = \|\hat \varphi - \varphi \|_{\Lsp}^2$.

\input{figs/network_architecture}

We note that this loss function is not invariant under phase shifts (i.e., multiplication by a scalar $c \in \mathbb C$ with $|c| = 1$). However, according to our numerical experiments, it often yields better results than the invariance-preserving loss function $\hat L(\hat \varphi , \varphi) = \| \hat \rho - \rho\|_{L^1(\mathcal D)}$.

\subsection{Acceleration strategy}\label{sec:strategy}

Our strategy crucially relies on the use of the normalization error of the network output as a criterion for its quality. Indeed, our neural network architecture explicitly does not implement a normalization layer. Thus, the output of the neural network $\tilde \varphi$ may have an arbitrary norm and we can use the normalization error
$$
e = |1 - \|\tilde \varphi \|_{\Lsp}|
$$
as an error indicator to estimate the quality of the neural network approximation. 

For every $n_{e}$-th iteration of the EARCG, if the energy norm of $g^{(k)}$ is in $[\epsilon_1^{\rm min}, \epsilon_1^{\rm max}]$, we generate a neural network prediction, which is not normalized in general. If its error indicator $e$ falls below the prescribed threshold $e_0$, we use the normalized output of the neural network, after normalization, as our new approximation and continue with EARCG until convergence.

Although the acceleration strategy could, in principle, be carried out multiple times during the algorithm, we limit ourselves to a single application. Additionally, if the energy norm of the gradient reaches the lower bound of the acceleration window $\epsilon_1^{\rm min}$ without any successful acceleration step, neural network acceleration is enforced. Although this occasionally leads to suboptimal results, it generally accelerates convergence based on our numerical experience. 

Figure~\ref{fig:enhanced-convergence} illustrates the capability of our approach by showcasing how a single well-timed acceleration, chosen according to our acceleration strategy, can overcome the slow realignment phase of pure EARCG outlined in the previous section in Figure~\ref{fig:slow-convergence}. 

\input{figs/heatmaps_convergence_enhanced}

\section{Numerical Experiments}\label{sec:numexp}

The numerical experiments use the EARCG implementation of \cite{PPS25}, which is based on the density functional toolkit \texttt{DFTK.jl} \cite{DFTKjcon}. Among others, it uses plane wave discretization, where the discretization level is determined by a kinetic energy cutoff $E_{\rm cut}$  of the plane waves, behaving as $E_{\rm cut} \sim 1/\sqrt{h}$ for grid size $h$. The source code of our experiments can be found at
\begin{center}
{
\vspace{0mm}
    \texttt{https://github.com/jonas-pueschel/NNacceleration4EARCG
    }
\vspace{0mm}
}
\end{center}
%
In what follows, we chose the parameters  $[\epsilon_1^{\rm min}, \epsilon_1^{\rm max}] = [10^{-4}, 10^{-1}]$, $\epsilon_2 = 10^{-8}$, $n_e = 5$ and $e_0 = 5 \times 10^{-3}$.

\subsection{Neural network training}
\begin{figure}
    \centering
    \includegraphics[width=\linewidth]{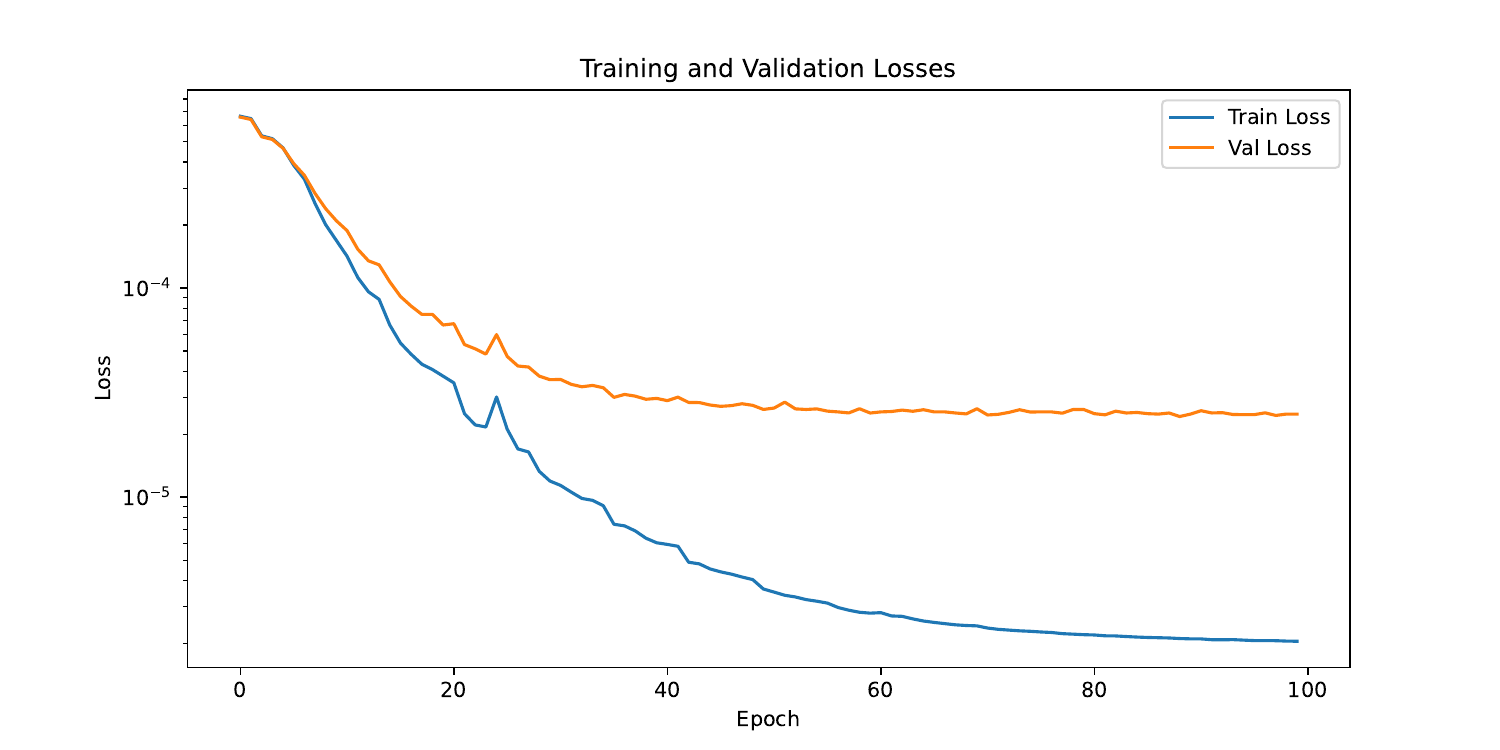}
    \caption{Training and validation loss for the UNet training.}
    \label{fig:losses}
\end{figure}

For the models considered, we fix the cutoff energy $E_{\rm cut} = 100$ and choose $a =20$ as well as $v_2 = 1$.
Always using new random initial data $\varphi_0 \in \SH$ for every data point, we generate two groups of data:
\begin{itemize}
    \item 1000 EARCG runs with parameters $\kappa, \omega, v_1$ sampled uniformly from
    $$\kappa \in [200, 1000],\quad \omega \in [0.8, 1.6], \quad v_1 \in [1,2]$$
    \item 1000 EARCG runs with parameters $\kappa, \omega, v_1$ sampled uniformly from
    $$\kappa \in [600, 1000],\quad \omega \in [1.2, 1.6], \quad v_1 \in [1,2]$$
\end{itemize}
The first set covers the whole parameter range we are interested in, while the second set focuses on the range resulting in more challenging systems.
From both sets, we sample $20$ data points for tolerances $\tilde \epsilon_j,\;j = 1, \dots, 20$ distributed log equidistant on $[\epsilon_1^{\rm min}, \epsilon_1^{\rm max}]$. More precisely, we get the formula
$$\tilde \epsilon_j = e^{ \left(1 - \tfrac{j-1}{19}\right)\ln(\epsilon_1^{\max}) + \tfrac{j-1}{19} \ln(\epsilon_1^{\min})}$$
The data points are tuples $( (\varphi_j, g_j), \varphi^\star )$ generated from the respective tolerances $\tilde \epsilon_j$ as described in Subsection \ref{sec:training_data}.
In total, this results in $40,000$ data points. 
The data is divided into $90$\% training and $10$\% evaluation data sets. Additionally, we apply classical data augmentation; for each drawn data point, perform a random horizontal and vertical flip (each with probability $p = 0.5$). Since the potentials are symmetric w.r.t. the horizontal and vertical axis through the origin, this allows the network to learn the symmetry. We train $100$ epochs of the network, and the losses are shown in Figure \ref{fig:losses}.

\subsection{Benchmarking}

Our analysis utilizes a test set of $500$ randomly generated initial states and parameters.

In order to show the efficacy of the neural network acceleration strategy laid out in Section \ref{sec:strategy}, we compare our approach to a single application of the network at a randomly chosen iteration. To this end, a tolerance $\epsilon_1$ is sampled log-uniformly from the interval $[\epsilon_1^{\min}, \epsilon_1^{\max}]$ and the network is applied once the EARCG reaches this tolerance. The neural network output is then used as input for a subsequent run of EARCG with tolerance $tol = \epsilon_2$. 

\input{figs/figure_hist_perc_its}

\input{figs/figure_hist_times}

\input{figs/figure_hist_rhos}

We utilize three quality measures to compare our strategy with random application.
Firstly, we calculate the difference in total steps between the classical and neural-accelerated algorithms. Secondly, we consider the percentage of iteration steps saved. 
Lastly, we use the relative improvement in density-error
\[
    \mathrm{impr}_\rho = \frac{\| \rho - \rho^*\|_{L^1(\Omega)} - \|\tilde \rho - \rho^*\|_{L^1(\Omega)}}{\| \rho - \rho^*\|_{L^1(\Omega)}},
\]
where $\varphi$ is the input to the network, $\tilde \varphi$ is its output, $\varphi^*$ is the reference solution (i.e., the EARCG converged state) and $\rho = | \varphi(\cdot)|^2$, $\tilde \rho = |\tilde \varphi(\cdot)|^2$, $\rho^* = |\varphi^*(\cdot)|^2$ their respective densities.
A value $\mathrm{impr}_\rho$ in $(0,1]$ indicates the neural network output is closer to the converged solution than the input, $0$ indicates no change in distance and negative values indicate the distance increased.
For all three measures, we additionally compute both the mean and median.

The histograms in Figures \ref{fig:figure_hist_perc_its}, \ref{fig:figure_hist_times} and \ref{fig:figure_hist_rhos} show the results. Both acceleration strategies reduce the number of iterations (absolute and relative), but our strategy achieves a more significant and reliable improvement, saving around $22\%$ of iterations on average compared to only minimal improvement for random application. Wall-clock time is reduced by $14.5\%$ on average. Acceleration reduces iterations in around $92\%$ of cases and wall time in $75\%$ of cases. The discrepancy between iterations and wall time arises because the neural network is applied every $n_e=5$ iterations until acceptance, incurring overhead. The algorithm is executed on a single CPU core in \texttt{julia}, and the neural network inference is also performed on the CPU. Leveraging a GPU could potentially reduce the overall wall time further.

Lastly, the error in density is reduced by around $46\%$ on average. 
The histograms in Figure \ref{fig:figure_hist_rhos} also show that sufficiently good density improvements usually guarantee convergence to the correct energy in the sense of the classical algorithm. Notably, cases where our acceleration strategy worsens the density approximation are extremely rare with $5.3\%$, compared to random application with $40\%$.

Among the $488$ cases where the classical method converged, random acceleration converged to the correct energy in $439$ cases. Our strategy succeeded in $463$ cases. Different final energies usually result from insufficient neural network acceleration or premature application, potentially leading to convergence to a different local minimum.

\subsection{Best and worst examples with respect to density-error}

To better understand the benefits and limitations of our approach, we examine three examples with the best and worst performance with respect to density-error.

We exclude cases where the classical iteration did not converge, i.e. the convergence criterion of residual norm smaller than $10^{-8}$ was not met within $30000$ iterations and consequently the resulting density of the classical method does not belong to a local minimizer. 
Figures \ref{fig:nn-best-rho-1} to \ref{fig:nn-best-rho-3} show the examples where density was improved the most, while Figures \ref{fig:nn-worst-rho-1} to \ref{fig:nn-worst-rho-3} show the ones where it was improved the least. 

The examples in Figures \ref{fig:nn-best-rho-1} to \ref{fig:nn-best-rho-3} demonstrate the ability of the neural network to identify the preferred rotation of the density, significantly improving the density in one step.  

In the examples shown in Figures \ref{fig:nn-worst-rho-1} the neural network in some sense is applied ``too late,'', leading to a deterioration of the density on paper. However, the number of iterations needed for convergence is still decreasing, accelerating actual convergence. Conversely, Figures \ref{fig:nn-worst-rho-2} and \ref{fig:nn-worst-rho-3} show cases where the neural network prediction is erroneous, resulting in an incorrect vortex configuration.

\input{figs/nn-best-rho-1}

\input{figs/nn-best-rho-2}

\input{figs/nn-best-rho-3}

\input{figs/nn-worst-rho-1}

\input{figs/nn-worst-rho-2}

\input{figs/nn-worst-rho-3}

\subsection{Best and worst examples with respect to iterations }

Next we discuss the three best and worst improvements in iteration count. We exclude examples, where the classical and enhanced methods converged to different local minima, since local convergence behavior may vary significantly, skewing the result. We note, however, that this excludes some cases where the acceleration yields subpar results.  
Figures \ref{fig:nn-best-it-1} to \ref{fig:nn-best-it-3} show cases with the largest relative reduction in iterations, while Figures \ref{fig:nn-worst-it-1} to \ref{fig:nn-worst-it-3} show cases with the largest increase.

In the best cases (Figures \ref{fig:nn-best-it-1}--\ref{fig:nn-best-it-3}), the neural network identifies the correct orientation of the vortices, thus reducing the number of iterations significantly. In all cases, the classical and enhanced algorithm eventually converged to the same local minimum; however, in the example from Figure \ref{fig:nn-best-it-2}, the classical iteration convergence is so slow that the energy is off by more than $10^{-8}$. 

In the worst cases (Figures \ref{fig:nn-worst-it-1}--\ref{fig:nn-worst-it-2}), the neural network fails to predict the vortex configuration correctly, leading to an increase in necessary iterations. In Figure \ref{fig:nn-worst-it-3}, the neural network predicts the vortex pattern correctly, but misaligns the orientation, whereas the classical method achieves correct orientation faster. We note the density does not deteriorate significantly in these cases; in two of the three, it even improves.

\input{figs/nn-best-it-1}

\input{figs/nn-best-it-2}

\input{figs/nn-best-it-3}

\input{figs/nn-worst-it-1}

\input{figs/nn-worst-it-2}

\input{figs/nn-worst-it-3}

\section{Conclusion and Open Problems}\label{sec:conclusion}

In this paper, we presented a neural network acceleration strategy for iterative Riemannian methods for energy minimization of the rotating Gross-Pitaevskii equation. Numerical experiments showed that using an error indicator based on the normalization error of the neural network enables saving $22\%$ of iterations and $14.5\%$ of wall time on average, resulting in significant acceleration.

Our approach may be further improved by allowing multiple neural network accelerations in a single EARCG run. It can also easily be extended to other descent methods if the required input state and descent direction are available, e.g. to different nonlinear eigenvalue problems like Hartree-Fock or Kohn-Sham. Regarding the neural network design, incorporating rotational invariance might improve performance.

\bibliographystyle{plain}
\bibliography{refs}
\end{document}

%% file: preamble.tex
\usepackage[toc,page]{appendix}
\usepackage{graphicx} 
\usepackage{url}
\usepackage{amsmath} 
\usepackage{amssymb}
\usepackage{bm, bbm}
\usepackage{polynom}
\usepackage[outline]{contour}
\usepackage{xcolor}
\usepackage{pgfplots}
\usepackage{mathtools}
\usepackage{stmaryrd}
\usepackage{diagbox}
\usepackage{subcaption}
\pgfplotsset{
  compat=1.13,
}
\usepackage{xfrac}  
\usepackage{enumitem} 
\usepackage{cleveref}
\usepackage{float}
\usepackage{comment}
\usepackage{mathrsfs}
\usepackage[automark]{scrlayer-scrpage}
\usepackage[ruled,vlined, algo2e]{algorithm2e}

\usepackage{tikz}
\usepackage{tikz-cd}
\usetikzlibrary{patterns}
\usepackage{scalerel,stackengine}

\usepackage{natbib}

\theoremstyle{plain}

\contourlength{0.1pt}
\contournumber{10}%
\definecolor{clr1}{RGB}{27,158,119}
\definecolor{clr2}{RGB}{217,95,2}
\definecolor{clr3}{RGB}{117,112,179}
\definecolor{clr4}{RGB}{231,41,138}
\definecolor{clr5}{RGB}{102,166,30}
\definecolor{clr6}{RGB}{230,171,2}
\definecolor{clr7}{RGB}{166,118,29}
\definecolor{clr8}{RGB}{102,102,102}

\newcommand\C{\mathbb C}

\newcommand{\A}{\mathcal{A}}
\newcommand{\E}{\mathcal{E}}

\newcommand{\R}{\mathcal{R}}
\newcommand{\T}{\mathcal{T}}

\newcommand{\IR}{\mathbb R}

\newcommand{\Tan}{\mathrm{T}}

\newcommand{\grad}{\mathrm{grad}}

\newcommand{\Drm}{\mathrm{D}\,}

\newcommand{\re}{\mbox{\rm Re}}

\usepackage[most]{tcolorbox}

\definecolor{myGreen2}{RGB}{114,175,30} 
\definecolor{myRed}{RGB}{180,50,50}  
\definecolor{myOrange}{RGB}{225,92,22} 

\definecolor{unia-purple}{RGB}{173, 0, 124}
\definecolor{light-gray}{rgb}{0.8889,0.8889,0.8889}
\definecolor{cpl1}{rgb}{0.8889,0.4356,0.2781}
\definecolor{cpl2}{rgb}{0.0,0.6056,0.9787}
\definecolor{cpl3}{rgb}{0.2422,0.6433,0.3044}

\newlength{\plotwidth}
\newlength{\plotheight}
\usetikzlibrary{arrows.meta}
\usetikzlibrary{fit}
\usetikzlibrary{calc}
\usetikzlibrary{arrows,shapes,trees, patterns, backgrounds}
\usepgfplotslibrary{patchplots}
\usepgfplotslibrary{fillbetween}
\pgfplotsset{%
    layers/standard/.define layer set={%
        background,axis background,axis grid,axis ticks,axis lines,axis tick labels,pre main,main,axis descriptions,axis foreground%
    }{
        grid style={/pgfplots/on layer=axis grid},%
        tick style={/pgfplots/on layer=axis ticks},%
        axis line style={/pgfplots/on layer=axis lines},%
        label style={/pgfplots/on layer=axis descriptions},%
        legend style={/pgfplots/on layer=axis descriptions},%
        title style={/pgfplots/on layer=axis descriptions},%
        colorbar style={/pgfplots/on layer=axis descriptions},%
        ticklabel style={/pgfplots/on layer=axis tick labels},%
        axis background@ style={/pgfplots/on layer=axis background},%
        3d box foreground style={/pgfplots/on layer=axis foreground},%
    },
}

\newcommand{\Hsp}{H^1}
\newcommand{\Hspi}{H^{-1}}

\newcommand{\Lsp}{L^2}

\newcommand{\SH}{\mathcal S}
\def\dx{\,\text{d}x}

\makeatletter
\newsavebox{\@brx}
\newcommand{\llangle}[1][]{\savebox{\@brx}{\(\m@th{#1\langle}\)}%
  \mathopen{\copy\@brx\mkern2mu\kern-0.9\wd\@brx\usebox{\@brx}}}
\newcommand{\rrangle}[1][]{\savebox{\@brx}{\(\m@th{#1\rangle}\)}%
  \mathclose{\copy\@brx\mkern2mu\kern-0.9\wd\@brx\usebox{\@brx}}}
\makeatother

\newcommand\calE{\mathcal E}

%% file: figs/heatmaps_convergence_classical.tex
\begin{figure}[H]
    \begin{tabular}{cccc}
        \includegraphics[width = 0.22\textwidth, trim={3.8cm 0.2cm 3.8cm 0.7cm},clip]{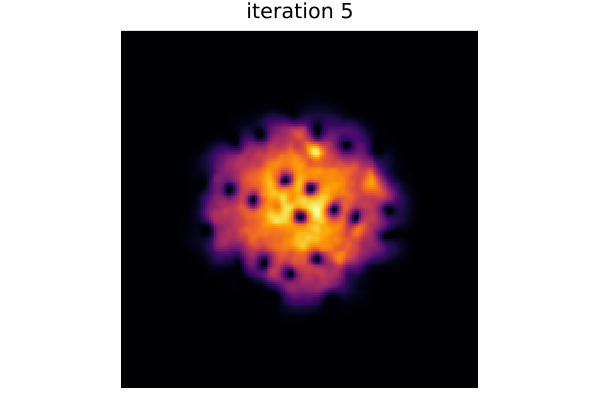} & 
        \includegraphics[width = 0.22\textwidth, trim={3.8cm 0.2cm 3.8cm 0.7cm},clip]{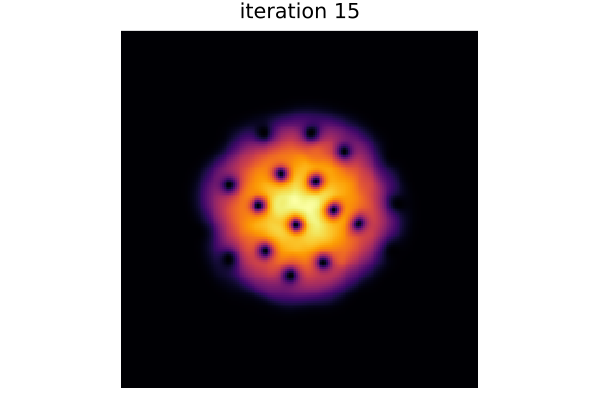} & 
        \includegraphics[width = 0.22\textwidth, trim={3.8cm 0.2cm 3.8cm 0.7cm},clip]{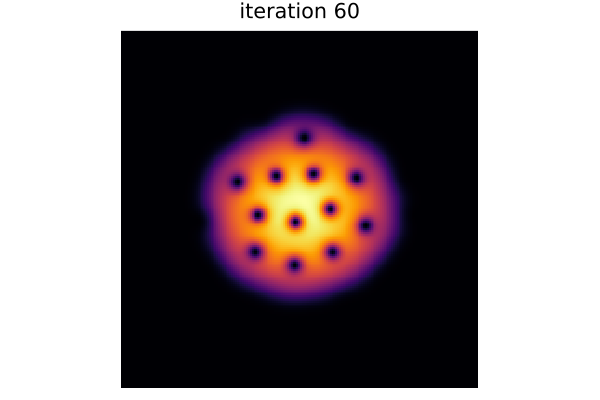} & 
        \includegraphics[width = 0.22\textwidth, trim={3.8cm 0.2cm 3.8cm 0.7cm},clip]{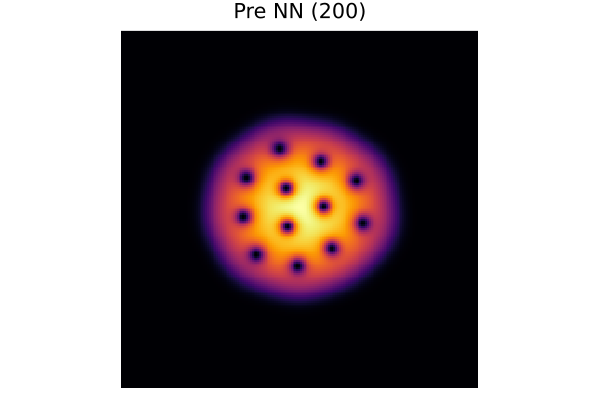} \\
        Iteration 5 & Iteration 15 & Iteration 60 & Iteration 200  \\[0.2cm]
        \includegraphics[width = 0.22\textwidth, trim={3.8cm 0.2cm 3.8cm 0.7cm},clip]{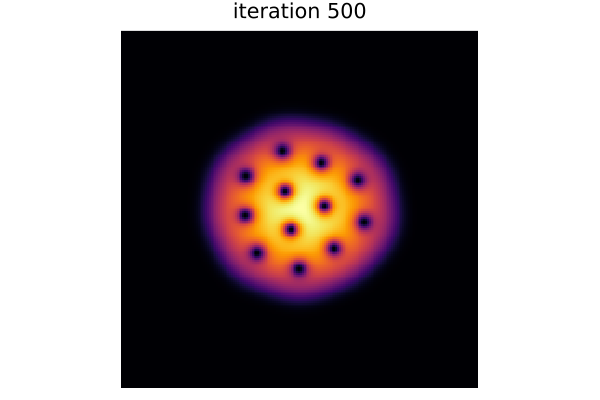} &
        \includegraphics[width = 0.22\textwidth, trim={3.8cm 0.2cm 3.8cm 0.7cm},clip]{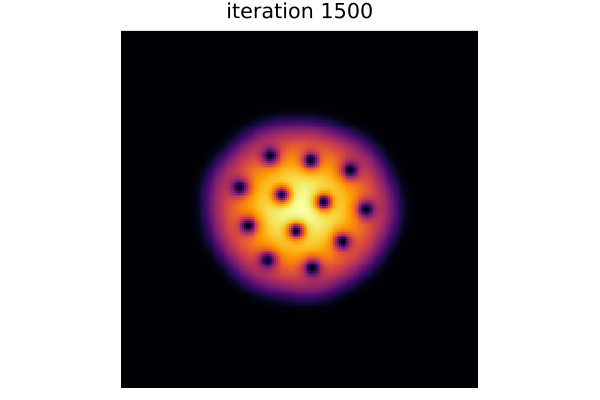} &
        \includegraphics[width = 0.22\textwidth, trim={3.8cm 0.2cm 3.8cm 0.7cm},clip]{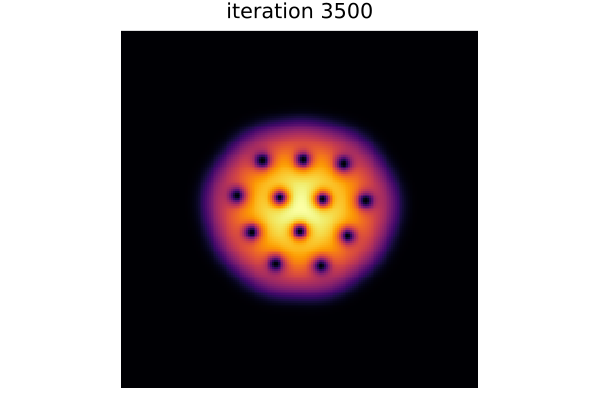} &
        \includegraphics[width = 0.22\textwidth, trim={3.8cm 0.2cm 3.8cm 0.7cm},clip]{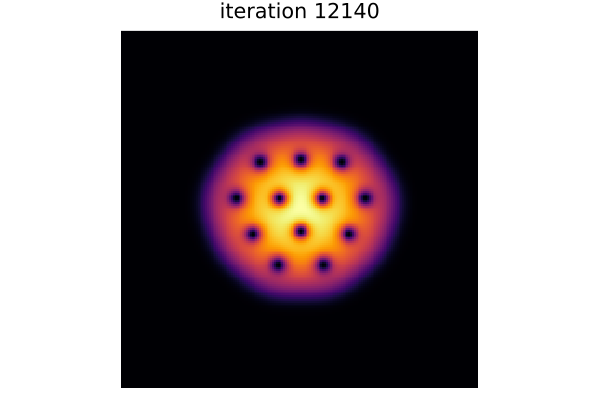} \\
        Iteration 500 & Iteration 1500 & Iteration 3500 & Iteration 12140
    \end{tabular}
    \caption{Density plots from EARCG iteration for parameters $a = 20, v = (1.1, 1), \omega = 1.4$ and $\kappa = 1000$. The final vortex configuration is reached after around $500$ iterations, convergence with the correct orientation is reached after $12140$ iterations.}
    \label{fig:slow-convergence}
\end{figure}

%% file: figs/network_architecture.tex
\pgfmathsetmacro{\xsep}{0.3}
\pgfmathsetmacro{\ysep}{0.4}
\pgfmathsetmacro{\hI}{1.9} 
\pgfmathsetmacro{\wI}{0.02} 
\pgfmathsetmacro{\wII}{0.08} 
\pgfmathsetmacro{\hII}{1.2} 
\pgfmathsetmacro{\wIII}{0.16} 
\pgfmathsetmacro{\hIII}{0.7} 
\pgfmathsetmacro{\wIV}{0.32} 
\pgfmathsetmacro{\hIV}{0.4} 
\pgfmathsetmacro{\wV}{0.64} 
\pgfmathsetmacro{\hV}{0.25} 
\pgfmathsetmacro{\wVI}{1.28} 
\pgfmathsetmacro{\hVI}{0.15} 
\pgfmathsetmacro{\shiftlabel}{0.25} 
\pgfmathsetmacro{\sizearrow}{0.3} 
\pgfmathsetmacro{\hseplegend}{0.23} 

\definecolor{mygrayfill}{RGB}{200,200,200}
\definecolor{mygraydraw}{RGB}{180,180,180}
\definecolor{mybluedashedfill}{RGB}{102,178,255}
\definecolor{mybluedasheddraw}{RGB}{51,153,255}
\definecolor{myyellowfill}{RGB}{255,255,204}
\definecolor{myyellowdraw}{RGB}{255,255,51}
\definecolor{myorangetext}{RGB}{255,153,51}

\definecolor{myvioletarrow}{RGB}{178,102,255}
\definecolor{mygreenarrow}{RGB}{0,205,0}
\definecolor{myredarrow}{RGB}{255,51,51}
\definecolor{mycyanarrow}{RGB}{0,204,204}
\definecolor{myorangearrow}{RGB}{255,128,0}
\definecolor{mypinkarrow}{RGB}{255,204,229}
\definecolor{mybluearrow}{RGB}{0,102,204}
\definecolor{mydashedarrow}{RGB}{255,153,51} 

\begin{figure}
\centering
\begin{tikzpicture}[
myrect/.style={
  rectangle,
  draw=mygraydraw,
  fill=mygrayfill,
  inner sep=0pt,
  fit=#1},
mydashedrect/.style={
  rectangle,
  draw=mybluedasheddraw,
  pattern=north west lines, 
  pattern color=mybluedashedfill, 
  dashed,
  inner sep=0pt, 
  fit=#1},
myinvisiblerect/.style={
  rectangle,
  draw=none,
  inner sep=0pt, 
  fit=#1},
myblockrect/.style={
  rectangle,
  draw=myyellowdraw,
  fill=myyellowfill,
  dashed,
  inner sep=0pt, 
  fit=#1},
upsampling/.style={
	->,
	mygreenarrow,
	thick},
bottleneck/.style={
	->,
	myvioletarrow,
	thick},
bottleneck_s2/.style={
	->,
	myredarrow,
	thick},
conv1x1/.style={
	->,
	mycyanarrow,
	thick},
conv3x3/.style={
	->,
	mybluearrow,
	thick},
conv7x7_s2/.style={
	->,
	mypinkarrow,
	thick},
maxpool2/.style={
	->,
	myredarrow,
	thick},
copy/.style={
	->,
	mydashedarrow, 
	thick, dashed},
labelfilters/.style 2 args={
	label={[black!60, shift={(#1,0.1)}]above:{\tiny #2}}}
][y=-1cm]

	\begin{scope}[y=-1cm]
		\coordinate (A_origin) at (0, 0);
		\coordinate (B_input) at ([shift={(\wI, \hI)}]A_origin);
		\node[myrect={(A_origin) (B_input)}, labelfilters={0}{4}] (input) {};

		\coordinate (A_b00u1) at  ([shift={(\wI + \xsep, 0)}]A_origin);
		\coordinate (B_b00u1) at ([shift={(\wII, \hI)}]A_b00u1);
		\node[myrect={(A_b00u1) (B_b00u1)}, labelfilters={0}{64}] (b00_u1) {};

        \coordinate (A_b00u2) at  ([shift={(\wII + \xsep, 0)}]A_b00u1);
		\coordinate (B_b00u2) at ([shift={(\wII, \hI)}]A_b00u2);
		\node[myrect={(A_b00u2) (B_b00u2)}, labelfilters={0}{64}] (b00_u2) {};
		
		\coordinate (A_b01u1) at ([shift={(0, \hI + \ysep)}]A_b00u2);
		\coordinate (B_b01u1) at ([shift={(\wII, \hII)}]A_b01u1);
		\node[myrect={(A_b01u1) (B_b01u1)}, labelfilters={-\shiftlabel}{64}] (b01_u1) {};

        \coordinate (A_b01u2) at ([shift={(\wII + \xsep, 0)}]A_b01u1);
		\coordinate (B_b01u2) at ([shift={(\wIII, \hII)}]A_b01u2);
		\node[myrect={(A_b01u2) (B_b01u2)}, labelfilters={0}{128}] (b01_u2) {};

        \coordinate (A_b01u3) at ([shift={(\wIII + \xsep, 0)}]A_b01u2);
		\coordinate (B_b01u3) at ([shift={(\wIII, \hII)}]A_b01u3);
		\node[myrect={(A_b01u3) (B_b01u3)}, labelfilters={0}{128}] (b01_u3) {};
		
		\coordinate (A_b1u0) at ([shift={(0, \hII + \ysep)}]A_b01u3); 
		\coordinate (B_b1u0) at ([shift={(\wIII, \hIII)}]A_b1u0); 
		\node[myrect={(A_b1u0) (B_b1u0)}, labelfilters={-\shiftlabel * 1.1}{128}] (b1_u0) {};
		
		\coordinate (A_b1u1) at ([shift={(\wIII + \xsep, 0)}]A_b1u0); 
		\coordinate (B_b1u1) at ([shift={(\wIV, \hIII)}]A_b1u1); 
		\node[myrect={(A_b1u1) (B_b1u1)}, labelfilters={0}{256}] (b1_u1) {};
		
		\coordinate (A_b1u2) at ([shift={(\wIV + \xsep, 0)}]A_b1u1); 
		\coordinate (B_b1u2) at ([shift={(\wIV, \hIII)}]A_b1u2); 
		\node[myrect={(A_b1u2) (B_b1u2)}, labelfilters={0}{256}] (b1_u2) {};
		
		\coordinate (A_b1u3) at ([shift={(0, \hIII + \ysep)}]A_b1u2); 
		\coordinate (B_b1u3) at ([shift={(\wIV, \hIV)}]A_b1u3); 
		\node[myrect={(A_b1u3) (B_b1u3)}, labelfilters={-\shiftlabel * 1.1}{256}] (b1_u3) {};
		
		\coordinate (A_b2u1) at ([shift={(\wIV + \xsep, 0)}]A_b1u3); 
		\coordinate (B_b2u1) at ([shift={(\wV, \hIV)}]A_b2u1); 
		\node[myrect={(A_b2u1) (B_b2u1)}, labelfilters={0}{512}] (b2_u1) {};
		
		\coordinate (A_b2u2) at ([shift={(\wIV + 2*\xsep, 0)}]A_b2u1); 
		\coordinate (B_b2u2) at ([shift={(\wV, \hIV)}]A_b2u2); 
		\node[myrect={(A_b2u2) (B_b2u2)}, labelfilters={0}{512}] (b2_u2) {};
		
		\coordinate (A_b2u4) at ([shift={(0, \hIV + \ysep)}]A_b2u2); 
		\coordinate (B_b2u4) at ([shift={(\wV, \hV)}]A_b2u4); 
		\node[myrect={(A_b2u4) (B_b2u4)}, labelfilters={-\shiftlabel * 1.1}{512}] (b2_u4) {};
		
		\coordinate (A_b3u1) at ([shift={(\wV + \xsep, 0)}]A_b2u4); 
		\coordinate (B_b3u1) at ([shift={(\wVI, \hV)}]A_b3u1); 
		\node[myrect={(A_b3u1) (B_b3u1)}, labelfilters={0}{1024}] (b3_u1) {};
		
		\coordinate (A_b3u2) at ([shift={(\wVI + \xsep, 0)}]A_b3u1); 
		\coordinate (B_b3u2) at ([shift={(\wVI, \hV)}]A_b3u2); 
		\node[myrect={(A_b3u2) (B_b3u2)}, labelfilters={\shiftlabel * 1.3}{1024}] (b3_u2) {};

        
		\coordinate (A_up2) at ([shift={(\wIV, -(\ysep + \hIV ))}]A_b3u2); 
		\coordinate (B_up2) at ([shift={(\wV, \hIV)}]A_up2); 
		\node[myrect={(A_up2) (B_up2)}, labelfilters={-\wIV}{1024}] (up2) {};
		\coordinate (A_cat2) at ([shift={(-\wV, 0)}]A_up2); 
		\coordinate (B_cat2) at ([shift={(\wV, \hIV)}]A_cat2); 
		\node[mydashedrect={(A_cat2) (B_cat2)}] (cat2) {};

        \coordinate (A_deconv2u1) at ([shift={(\wV + \xsep, 0)}]A_up2); 
		\coordinate (B_deconv2u1) at ([shift={(\wV, \hIV)}]A_deconv2u1); 
		\node[myrect={(A_deconv2u1) (B_deconv2u1)}, labelfilters={0}{512}] (deconv2_1) {};
        
        \coordinate (A_deconv2u2) at ([shift={(\wV + \xsep, 0)}]A_deconv2u1); 
		\coordinate (B_deconv2u2) at ([shift={(\wV, \hIV)}]A_deconv2u2); 
		\node[myrect={(A_deconv2u2) (B_deconv2u2)}, labelfilters={\shiftlabel}{512}] (deconv2_2) {};
		
		\coordinate (A_up3) at ([shift={(\wIII, -(\ysep + \hIII))}]A_deconv2u2); 
		\coordinate (B_up3) at ([shift={(\wIV, \hIII)}]A_up3); 
		\node[myrect={(A_up3) (B_up3)}, labelfilters={-\wIII}{512}] (up3) {};
		\coordinate (A_cat3) at ([shift={(-\wIV, 0)}]A_up3); 
		\coordinate (B_cat3) at ([shift={(\wIV, \hIII)}]A_cat3); 
		\node[mydashedrect={(A_cat3) (B_cat3)}] (cat3) {};

        \coordinate (A_deconv3u1) at ([shift={(\wIV + \xsep, 0)}]A_up3); 
		\coordinate (B_deconv3u1) at ([shift={(\wIV, \hIII)}]A_deconv3u1); 
		\node[myrect={(A_deconv3u1) (B_deconv3u1)}, labelfilters={0}{256}] (deconv3_1) {};
        
        \coordinate (A_deconv3u2) at ([shift={(\wIV + \xsep, 0)}]A_deconv3u1); 
		\coordinate (B_deconv3u2) at ([shift={(\wIV, \hIII)}]A_deconv3u2); 
		\node[myrect={(A_deconv3u2) (B_deconv3u2)}, labelfilters={\shiftlabel}{256}] (deconv3_2) {};

		\coordinate (A_up4) at ([shift={(\wII, -(\ysep + \hII))}]A_deconv3u2); 
		\coordinate (B_up4) at ([shift={(\wIII, \hII)}]A_up4); 
		\node[myrect={(A_up4) (B_up4)}, labelfilters={-\wII}{256}] (up4) {};
		\coordinate (A_cat4) at ([shift={(-\wIII, 0)}]A_up4); 
		\coordinate (B_cat4) at ([shift={(\wIII, \hII)}]A_cat4); 
		\node[mydashedrect={(A_cat4) (B_cat4)}] (cat4) {};

        \coordinate (A_deconv4u1) at ([shift={(\wIII + \xsep, 0)}]A_up4); 
		\coordinate (B_deconv4u1) at ([shift={(\wIII, \hII)}]A_deconv4u1); 
		\node[myrect={(A_deconv4u1) (B_deconv4u1)}, labelfilters={0}{128}] (deconv4_1) {};
        
        \coordinate (A_deconv4u2) at ([shift={(\wIII + \xsep, 0)}]A_deconv4u1); 
		\coordinate (B_deconv4u2) at ([shift={(\wIII, \hII)}]A_deconv4u2); 
		\node[myrect={(A_deconv4u2) (B_deconv4u2)}, labelfilters={\shiftlabel}{128}] (deconv4_2) {};
		
		\coordinate (A_up5) at ([shift={(0.5 * \wII, -(\ysep + \hI))}]A_deconv4u2); 
		\coordinate (B_up5) at ([shift={(\wII, \hI)}]A_up5); 
		\node[myrect={(A_up5) (B_up5)}, labelfilters={-0.5*\wII}{128}] (up5) {};
		\coordinate (A_cat5) at ([shift={(-\wII, 0)}]A_up5); 
		\coordinate (B_cat5) at ([shift={(\wII, \hI)}]A_cat5); 
		\node[mydashedrect={(A_cat5) (B_cat5)}] (cat5) {};

        \coordinate (A_deconv5u1) at ([shift={(\wII + \xsep, 0)}]A_up5); 
		\coordinate (B_deconv5u1) at ([shift={(\wII, \hI)}]A_deconv5u1); 
		\node[myrect={(A_deconv5u1) (B_deconv5u1)}, labelfilters={0}{64}] (deconv5_1) {};
        
        \coordinate (A_deconv5u2) at ([shift={(\wII + \xsep, 0)}]A_deconv5u1); 
		\coordinate (B_deconv5u2) at ([shift={(\wII, \hI)}]A_deconv5u2); 
		\node[myrect={(A_deconv5u2) (B_deconv5u2)}, labelfilters={0}{64}] (deconv5_2) {};
		
		\coordinate (A_out) at ([shift={(\wII + \xsep, 0)}]A_deconv5u2); 
		\coordinate (B_out) at ([shift={(\wI, \hI)}]A_out); 
		\node[myrect={(A_out) (B_out)}, labelfilters={0}{2}] (output) {};
			
		\draw[conv3x3] (input) -- (b00_u1);
        \draw[conv3x3] (b00_u1) -- (b00_u2);
		\draw[maxpool2] (b00_u2) -- (b01_u1);
		\draw[conv3x3] (b01_u1) -- (b01_u2);
        \draw[conv3x3] (b01_u2) -- (b01_u3);
        \draw[maxpool2] (b01_u3) -- (b1_u0);
		
		\draw[conv3x3] (b1_u0) -- (b1_u1);
        \draw[conv3x3] (b1_u1) -- (b1_u2);
		\draw[maxpool2] (b1_u2) -- (b1_u3);

        \draw[conv3x3] (b1_u3) -- (b2_u1);
        \draw[conv3x3] (b2_u1) -- (b2_u2);
		\draw[maxpool2] (b2_u2) -- (b2_u4);
		\draw[conv3x3] (b2_u4) -- (b3_u1);
		
		\draw[conv3x3] (b3_u1) -- (b3_u2);
		
        \draw[upsampling] (b3_u2) -- (up2);

        \draw[conv3x3] (up2) -- (deconv2_1);
        \draw[conv3x3] (deconv2_1) -- (deconv2_2);
        
		\draw[upsampling] (deconv2_2) -- (up3);
		\draw[conv3x3] (up3) -- (deconv3_1);
        \draw[conv3x3] (deconv3_1) -- (deconv3_2);
		
		\draw[upsampling] (deconv3_2) -- (up4);
		\draw[conv3x3] (up4) -- (deconv4_1);
        \draw[conv3x3] (deconv4_1) -- (deconv4_2);
		
		\draw[upsampling] (deconv4_2) -- (up5);
		\draw[conv3x3] (up5) -- (deconv5_1);
        \draw[conv3x3] (deconv5_1) -- (deconv5_2);
		\draw[copy] (b2_u2) -- (cat2);
		\draw[copy] (b1_u2) -- (cat3);
		\draw[copy] (b01_u3) -- (cat4);
		\draw[copy] (b00_u2) -- (cat5);
		
		\draw[conv1x1] (deconv5_2) -- (output);

		\coordinate (A_s) at ([shift={(0.5*\xsep, 0)}]A_out); 
		\coordinate (B_s) at ([shift={(0, \hI)}]A_s); 
		\node [myinvisiblerect={(A_s) (B_s)}, label={[gray]right:{\tiny $n$}}] {};
		
		\coordinate (A_s2) at ([shift={(0, \hI + \ysep)}]A_s); 
		\coordinate (B_s2) at ([shift={(0, \hII)}]A_s2); 
		\node [myinvisiblerect={(A_s2) (B_s2)}, label={[gray]right:{\tiny $n/2$}}] {};
		
		\coordinate (A_s4) at ([shift={(0, \hII + \ysep)}]A_s2); 
		\coordinate (B_s4) at ([shift={(0, \hIII)}]A_s4); 
		\node [myinvisiblerect={(A_s4) (B_s4)}, label={[gray]right:{\tiny $n/4$}}] {};
		
		\coordinate (A_s8) at ([shift={(0, \hIII + \ysep)}]A_s4); 
		\coordinate (B_s8) at ([shift={(0, \hIV)}]A_s8); 
		\node [myinvisiblerect={(A_s8) (B_s8)}, label={[gray]right:{\tiny $n/8$}}] {};
		
		\coordinate (A_s16) at ([shift={(0, \hIV + \ysep)}]A_s8); 
		\coordinate (B_s16) at ([shift={(0, \hV)}]A_s16); 
		\node [myinvisiblerect={(A_s16) (B_s16)}, label={[gray]right:{\tiny $n/16$}}] {};
		
	\end{scope}

\coordinate (O_legend) at (0, -4.9);
\draw [maxpool2] ([shift={(0, -\hseplegend)}]O_legend) to ([shift={(+\sizearrow, -\hseplegend)}]O_legend) node[right] {\tiny max pool 2x2};
\draw [conv1x1] ([shift={(0, -2*\hseplegend)}]O_legend) to ([shift={(+\sizearrow, -2*\hseplegend)}]O_legend) node[right] {\tiny conv 1x1};
\draw [upsampling] ([shift={(0, -3*\hseplegend)}]O_legend) to ([shift={(+\sizearrow, -3*\hseplegend)}]O_legend) node[right] {\tiny up-conv 2x2};
\draw [conv3x3] ([shift={(0, -4*\hseplegend)}]O_legend) to ([shift={(+\sizearrow, -4*\hseplegend)}]O_legend) node[right] {\tiny conv 3x3};
\draw [copy] ([shift={(0, -5*\hseplegend)}]O_legend) to ([shift={(+\sizearrow, -5*\hseplegend)}]O_legend) node[right] {\tiny copy};

\end{tikzpicture}
\caption{Network architecture of the U-Net with a $4$-channel input and a $2$-channel output. We use a standard U-Net architecture with five blocks in the contracting and expanding paths respectively. Every block consists of two $3\times3$ padded convolution layers using ReLU activation functions. Contracting paths use $2\times2$ max-pooling (halving length and width of the tensor), while expanding paths use $2\times 2$ padded up-convolutions (doubling the length and width of the tensor). Skip connections are used send information directly from the contracting to the expanding path. For the output, an additional $1\times 1$ convolution layer is used.} 
\label{tikz:network}
\end{figure}

%% file: figs/heatmaps_convergence_enhanced.tex
\begin{figure}[H]
    \begin{tabular}{cccc}
        \includegraphics[width = 0.22\textwidth, trim={3.8cm 0.2cm 3.8cm 0.7cm},clip]{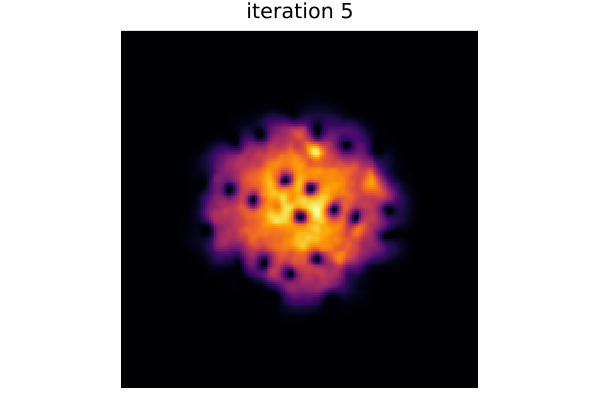} & 
        \includegraphics[width = 0.22\textwidth, trim={3.8cm 0.2cm 3.8cm 0.7cm},clip]{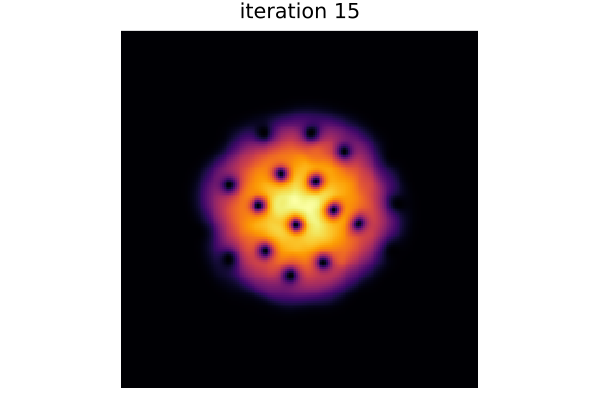} & 
        \includegraphics[width = 0.22\textwidth, trim={3.8cm 0.2cm 3.8cm 0.7cm},clip]{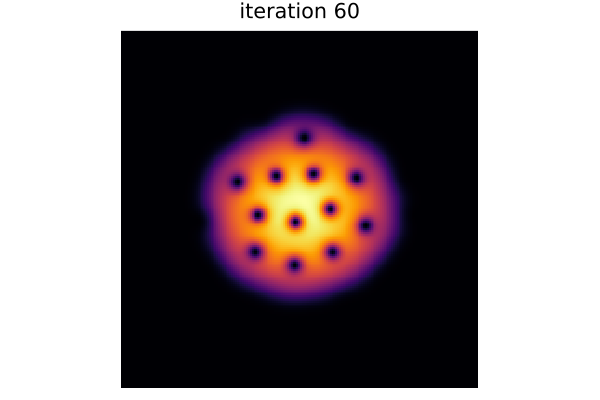} & 
        \includegraphics[width = 0.22\textwidth, trim={3.8cm 0.2cm 3.8cm 0.7cm},clip]{figs/heatmaps_convergence_bigbox/hmp-enh-pre.png} \\
        Iteration 5 & Iteration 15 & Iteration 60 & Pre NN (Iter. 200)  \\[0.2cm]
        \includegraphics[width = 0.22\textwidth, trim={3.8cm 0.2cm 3.8cm 0.7cm},clip]{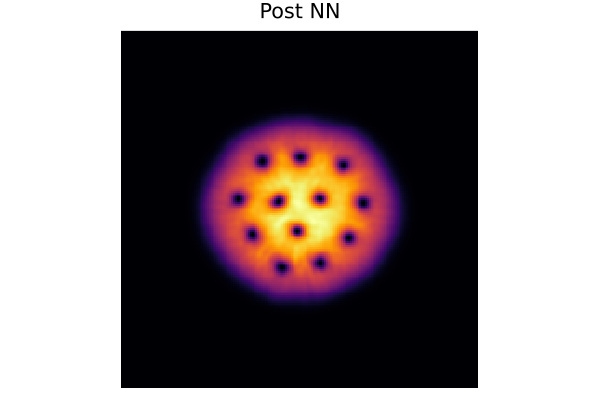} &
        \includegraphics[width = 0.22\textwidth, trim={3.8cm 0.2cm 3.8cm 0.7cm},clip]{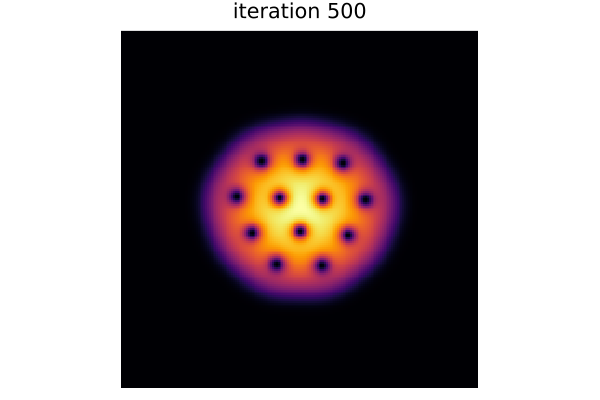} &
        \includegraphics[width = 0.22\textwidth, trim={3.8cm 0.2cm 3.8cm 0.7cm},clip]{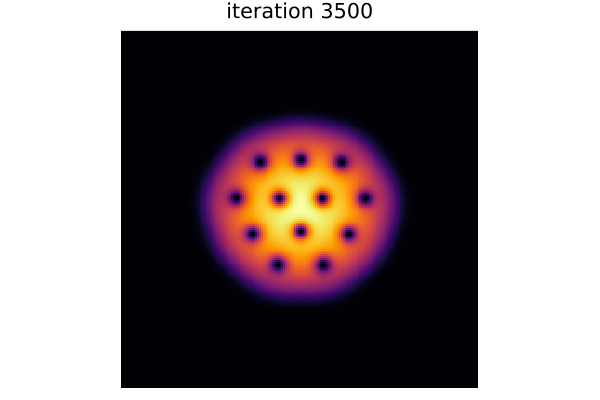} &
        \includegraphics[width = 0.22\textwidth, trim={3.8cm 0.2cm 3.8cm 0.7cm},clip]{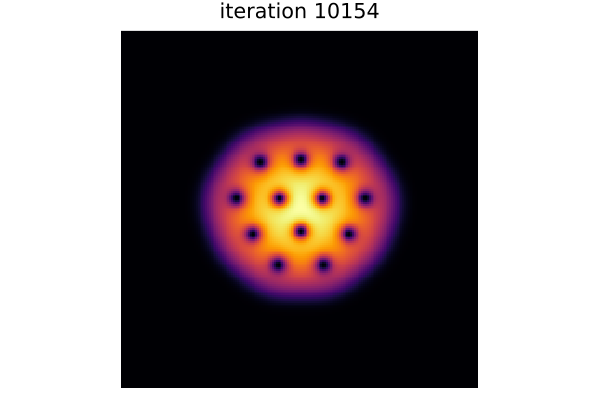} \\
        Post NN (Iter. 200)  & Iteration 500 & Iteration 3500 & Iteration 10152
    \end{tabular}
    \caption{Density plots from the neural network accelerated EARCG iteration for parameters $a = 20, v = (1.1, 1), \omega = 1.4$ and $\kappa = 1000$. Convergence is reached after $10152$ iterations, saving almost $2000$ iterations.}
    \label{fig:enhanced-convergence}
\end{figure}

%% file: figs/figure_hist_perc_its.tex
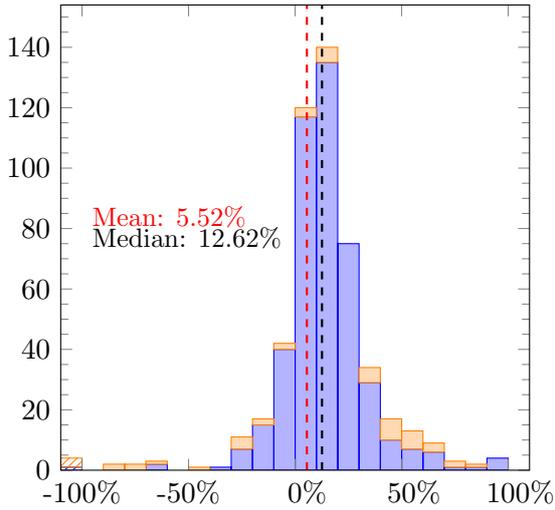
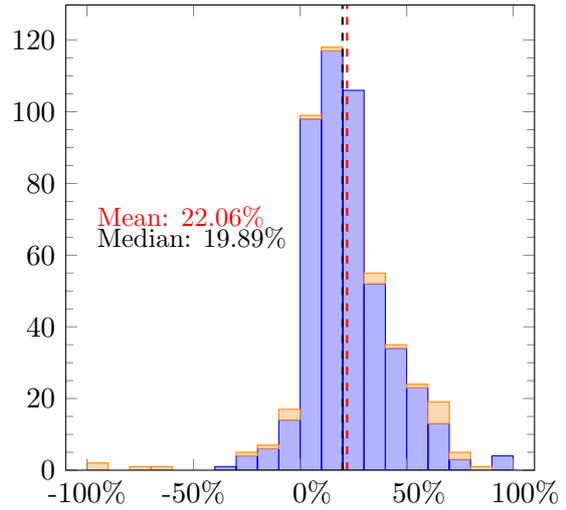
\begin{figure}
\centering
\begin{subfigure}{.5\textwidth}
  \centering
    \captionsetup{width=.9\textwidth}
  \input{./figs/histograms/hist_perc_its.tex}
  \label{fig:sub1}
    \subcaption{Percentage of iterations saved by neural network enhancement of the algorithm when neural network was applied at a randomly chosen iteration. 
    Improvement was achieved in 82.4\% of cases.
    }
\end{subfigure}%
\begin{subfigure}{.5\textwidth}
  \centering
  \captionsetup{width=.9\textwidth}
  \input{./figs/histograms/hist_perc_its_enhanced.tex}  
  \label{fig:sub2}
  \subcaption{Percentage of iterations saved by neural network enhancement of the algorithm when neural network was applied via the previously described acceleration strategy.
  Improvement was achieved in 91.8\% of cases.
  }
\end{subfigure}
    \caption{Histograms of the percentage of iterations saved by the neural network-accelerated algorithms, evaluated on all 500 test cases. 
    In blue are the cases, where the enhanced and classical algorithm converged to the same local minimum, in orange the cases of different minima. 
    The dashed bars denote the cases with increase bigger than 100\%.
    }
\label{fig:figure_hist_perc_its}
\end{figure}

%% file: figs/histograms/hist_perc_its.tex
\begin{tikzpicture}
\begin{axis}[
    width = \textwidth,
    height = \textwidth,
    ymin=0, ymax=154.0,
    xmin = -110, xmax = 110,
    minor y tick num = 3,
    xtick = {-100, -50, 0, 50, 100},
    xticklabels = {{-100\%, -50\%, \phantom{\%}0\%, \phantom{\%}50\%, \phantom{\%}100\%}},
    area style,
    ]
\addplot+[ybar stacked, bar width=10, color = blue!30, draw = blue] coordinates { 
(-105, 0)
(-95,0.0)
(-85,0.0)
(-75,0.0)
(-65,2.0)
(-55,0.0)
(-45,0.0)
(-35,1.0)
(-25,7.0)
(-15,15.0)
(-5,40.0)
(5,117.0)
(15,135.0)
(25,75.0)
(35,29.0)
(45,10.0)
(55,7.0)
(65,6.0)
(75,1.0)
(85,1.0)
(95,4.0)
(105, 0)
};
\addplot+[ybar stacked, bar width=10, color = blue!30, draw = blue, pattern = north east lines, pattern color = blue] coordinates { 
(-105, 1.0)
(-95,0)
(-85,0)
(-75,0)
(-65,0)
(-55,0)
(-45,0)
(-35,0)
(-25,0)
(-15,0)
(-5,0)
(5,0)
(15,0)
(25,0)
(35,0)
(45,0)
(55,0)
(65,0)
(75,0)
(85,0)
(95,0)
(105, 0.0)
};
\addplot+[ybar stacked, bar width=10, color = orange!30, draw = orange] coordinates { 
(-105, 0)
(-95,0.0)
(-85,2.0)
(-75,2.0)
(-65,1.0)
(-55,0.0)
(-45,1.0)
(-35,0.0)
(-25,4.0)
(-15,2.0)
(-5,2.0)
(5,3.0)
(15,5.0)
(25,0.0)
(35,5.0)
(45,7.0)
(55,6.0)
(65,3.0)
(75,2.0)
(85,1.0)
(95,0.0)
(105, 0)
};
\addplot+[ybar stacked, bar width=10, color = orange!30, draw = orange, pattern = north east lines, pattern color = orange] coordinates { 
(-105, 3.0)
(-95,0)
(-85,0)
(-75,0)
(-65,0)
(-55,0)
(-45,0)
(-35,0)
(-25,0)
(-15,0)
(-5,0)
(5,0)
(15,0)
(25,0)
(35,0)
(45,0)
(55,0)
(65,0)
(75,0)
(85,0)
(95,0)
(105, 0.0)
};

\addplot[mark=none, red ,dashed, thick] coordinates {(5.519153086975221,0) (5.519153086975221,280.0)};
\addplot [red,, nodes near coords={\small Mean: 5.52\%},every node near coord/.style={anchor=180}] coordinates {( -100, 84.0)};
\addplot[mark=none, black ,dashed ,thick] coordinates {(12.617614042249208,0) (12.617614042249208,280.0)};
\addplot [black,, nodes near coords={\small Median: 12.62\%},every node near coord/.style={anchor=180}] coordinates {(-100, 77.0)};
\end{axis}
\end{tikzpicture}

%% file: figs/histograms/hist_perc_its_enhanced.tex
\begin{tikzpicture}
\begin{axis}[
    width = \textwidth,
    height = \textwidth,
    ymin=0, ymax=129.8,
    xmin = -110, xmax = 110,
    minor y tick num = 3,
    xtick = {-100, -50, 0, 50, 100},
    xticklabels = {{-100\%, -50\%, \phantom{\%}0\%, \phantom{\%}50\%, \phantom{\%}100\%}},
    area style,
    ]
\addplot+[ybar stacked, bar width=10, color = blue!30, draw = blue] coordinates { 
(-105, 0)
(-95,0.0)
(-85,0.0)
(-75,0.0)
(-65,0.0)
(-55,0.0)
(-45,0.0)
(-35,1.0)
(-25,4.0)
(-15,6.0)
(-5,14.0)
(5,98.0)
(15,117.0)
(25,106.0)
(35,52.0)
(45,34.0)
(55,23.0)
(65,13.0)
(75,3.0)
(85,0.0)
(95,4.0)
(105, 0)
};
\addplot+[ybar stacked, bar width=10, color = blue!30, draw = blue, pattern = north east lines, pattern color = blue] coordinates { 
(-105, 0.0)
(-95,0)
(-85,0)
(-75,0)
(-65,0)
(-55,0)
(-45,0)
(-35,0)
(-25,0)
(-15,0)
(-5,0)
(5,0)
(15,0)
(25,0)
(35,0)
(45,0)
(55,0)
(65,0)
(75,0)
(85,0)
(95,0)
(105, 0.0)
};
\addplot+[ybar stacked, bar width=10, color = orange!30, draw = orange] coordinates { 
(-105, 0)
(-95,2.0)
(-85,0.0)
(-75,1.0)
(-65,1.0)
(-55,0.0)
(-45,0.0)
(-35,0.0)
(-25,1.0)
(-15,1.0)
(-5,3.0)
(5,1.0)
(15,1.0)
(25,0.0)
(35,3.0)
(45,1.0)
(55,1.0)
(65,6.0)
(75,2.0)
(85,1.0)
(95,0.0)
(105, 0)
};
\addplot+[ybar stacked, bar width=10, color = orange!30, draw = orange, pattern = north east lines, pattern color = orange] coordinates { 
(-105, 0.0)
(-95,0)
(-85,0)
(-75,0)
(-65,0)
(-55,0)
(-45,0)
(-35,0)
(-25,0)
(-15,0)
(-5,0)
(5,0)
(15,0)
(25,0)
(35,0)
(45,0)
(55,0)
(65,0)
(75,0)
(85,0)
(95,0)
(105, 0.0)
};

\addplot[mark=none, red ,dashed, thick] coordinates {(22.064948870533605,0) (22.064948870533605,236.0)};
\addplot [red,, nodes near coords={\small Mean: 22.06\%},every node near coord/.style={anchor=180}] coordinates {( -100, 70.8)};
\addplot[mark=none, black ,dashed ,thick] coordinates {(19.888584569272908,0) (19.888584569272908,236.0)};
\addplot [black,, nodes near coords={\small Median: 19.89\%},every node near coord/.style={anchor=180}] coordinates {(-100, 64.9)};
\end{axis}
\end{tikzpicture}

%% file: figs/figure_hist_times.tex
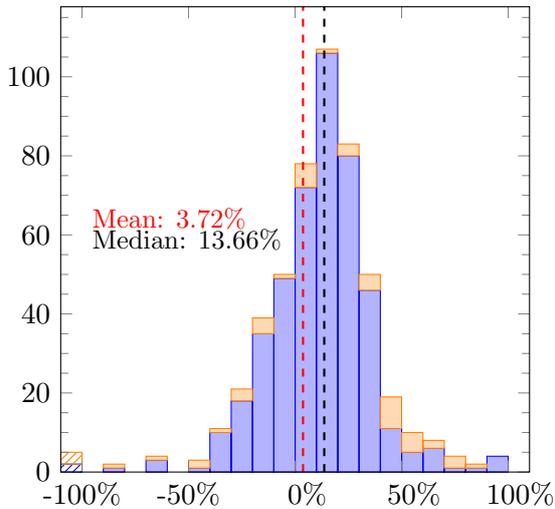
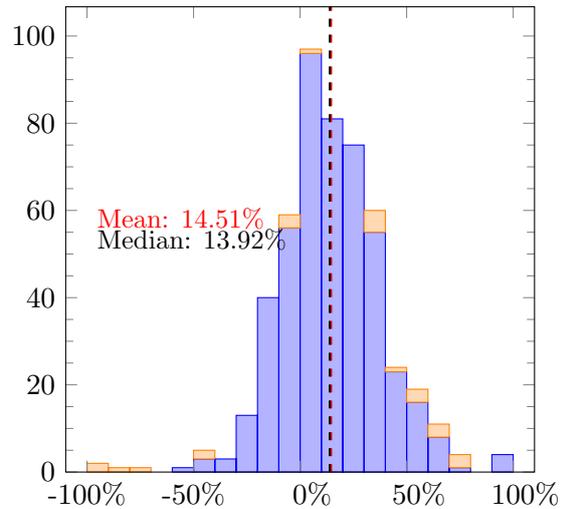
\begin{figure}
\centering
\begin{subfigure}{.5\textwidth}
  \centering
    \captionsetup{width=.9\textwidth}
  \input{./figs/histograms/hist_times.tex}
  \label{fig:sub1}
    \subcaption{Reduction in wall time by neural network enhancement of the algorithm when neural network was applied at a randomly chosen iteration. 
    Improvement was achieved in 73.0\% of cases.
    }
\end{subfigure}%
\begin{subfigure}{.5\textwidth}
  \centering
  \captionsetup{width=.9\textwidth}
  \input{./figs/histograms/hist_times_enhanced.tex}  
  \label{fig:sub2}
  \subcaption{Reduction in wall time by neural network enhancement of the algorithm when neural network was applied via the previously described acceleration strategy.
  Improvement was achieved in 75.0\% of cases.
  }
\end{subfigure}
    \caption{Histograms of the reduction in wall time by the neural network-accelerated algorithms, evaluated on all 500 test cases. 
    In blue are the cases, where the enhanced and classical algorithm converged to the same local minimum, in orange the cases of different minima. 
    The dashed bars denote the cases with increase bigger than 100\%.
    }
\label{fig:figure_hist_times}
\end{figure}

%% file: figs/histograms/hist_times.tex
\begin{tikzpicture}
\begin{axis}[
    width = \textwidth,
    height = \textwidth,
    ymin=0, ymax=117.7,
    xmin = -110, xmax = 110,
    minor y tick num = 3,
    xtick = {-100, -50, 0, 50, 100},
    xticklabels = {{-100\%, -50\%, \phantom{\%}0\%, \phantom{\%}50\%, \phantom{\%}100\%}},
    area style,
    ]
\addplot+[ybar stacked, bar width=10, color = blue!30, draw = blue] coordinates { 
(-105, 0)
(-95,0.0)
(-85,1.0)
(-75,0.0)
(-65,3.0)
(-55,0.0)
(-45,1.0)
(-35,10.0)
(-25,18.0)
(-15,35.0)
(-5,49.0)
(5,72.0)
(15,106.0)
(25,80.0)
(35,46.0)
(45,11.0)
(55,5.0)
(65,6.0)
(75,1.0)
(85,1.0)
(95,4.0)
(105, 0)
};
\addplot+[ybar stacked, bar width=10, color = blue!30, draw = blue, pattern = north east lines, pattern color = blue] coordinates { 
(-105, 2.0)
(-95,0)
(-85,0)
(-75,0)
(-65,0)
(-55,0)
(-45,0)
(-35,0)
(-25,0)
(-15,0)
(-5,0)
(5,0)
(15,0)
(25,0)
(35,0)
(45,0)
(55,0)
(65,0)
(75,0)
(85,0)
(95,0)
(105, 0.0)
};
\addplot+[ybar stacked, bar width=10, color = orange!30, draw = orange] coordinates { 
(-105, 0)
(-95,0.0)
(-85,1.0)
(-75,0.0)
(-65,1.0)
(-55,0.0)
(-45,2.0)
(-35,1.0)
(-25,3.0)
(-15,4.0)
(-5,1.0)
(5,6.0)
(15,1.0)
(25,3.0)
(35,4.0)
(45,8.0)
(55,5.0)
(65,2.0)
(75,3.0)
(85,1.0)
(95,0.0)
(105, 0)
};
\addplot+[ybar stacked, bar width=10, color = orange!30, draw = orange, pattern = north east lines, pattern color = orange] coordinates { 
(-105, 3.0)
(-95,0)
(-85,0)
(-75,0)
(-65,0)
(-55,0)
(-45,0)
(-35,0)
(-25,0)
(-15,0)
(-5,0)
(5,0)
(15,0)
(25,0)
(35,0)
(45,0)
(55,0)
(65,0)
(75,0)
(85,0)
(95,0)
(105, 0.0)
};

\addplot[mark=none, red ,dashed, thick] coordinates {(3.715359717433032,0) (3.715359717433032,214.0)};
\addplot [red,, nodes near coords={\small Mean: 3.72\%},every node near coord/.style={anchor=180}] coordinates {( -100, 64.2)};
\addplot[mark=none, black ,dashed ,thick] coordinates {(13.664313889626767,0) (13.664313889626767,214.0)};
\addplot [black,, nodes near coords={\small Median: 13.66\%},every node near coord/.style={anchor=180}] coordinates {(-100, 58.85)};
\end{axis}
\end{tikzpicture}

%% file: figs/histograms/hist_times_enhanced.tex
\begin{tikzpicture}
\begin{axis}[
    width = \textwidth,
    height = \textwidth,
    ymin=0, ymax=106.7,
    xmin = -110, xmax = 110,
    minor y tick num = 3,
    xtick = {-100, -50, 0, 50, 100},
    xticklabels = {{-100\%, -50\%, \phantom{\%}0\%, \phantom{\%}50\%, \phantom{\%}100\%}},
    area style,
    ]
\addplot+[ybar stacked, bar width=10, color = blue!30, draw = blue] coordinates { 
(-105, 0)
(-95,0.0)
(-85,0.0)
(-75,0.0)
(-65,0.0)
(-55,1.0)
(-45,3.0)
(-35,3.0)
(-25,13.0)
(-15,40.0)
(-5,56.0)
(5,96.0)
(15,81.0)
(25,75.0)
(35,55.0)
(45,23.0)
(55,16.0)
(65,8.0)
(75,1.0)
(85,0.0)
(95,4.0)
(105, 0)
};
\addplot+[ybar stacked, bar width=10, color = blue!30, draw = blue, pattern = north east lines, pattern color = blue] coordinates { 
(-105, 0.0)
(-95,0)
(-85,0)
(-75,0)
(-65,0)
(-55,0)
(-45,0)
(-35,0)
(-25,0)
(-15,0)
(-5,0)
(5,0)
(15,0)
(25,0)
(35,0)
(45,0)
(55,0)
(65,0)
(75,0)
(85,0)
(95,0)
(105, 0.0)
};
\addplot+[ybar stacked, bar width=10, color = orange!30, draw = orange] coordinates { 
(-105, 0)
(-95,2.0)
(-85,1.0)
(-75,1.0)
(-65,0.0)
(-55,0.0)
(-45,2.0)
(-35,0.0)
(-25,0.0)
(-15,0.0)
(-5,3.0)
(5,1.0)
(15,0.0)
(25,0.0)
(35,5.0)
(45,1.0)
(55,3.0)
(65,3.0)
(75,3.0)
(85,0.0)
(95,0.0)
(105, 0)
};
\addplot+[ybar stacked, bar width=10, color = orange!30, draw = orange, pattern = north east lines, pattern color = orange] coordinates { 
(-105, 0.0)
(-95,0)
(-85,0)
(-75,0)
(-65,0)
(-55,0)
(-45,0)
(-35,0)
(-25,0)
(-15,0)
(-5,0)
(5,0)
(15,0)
(25,0)
(35,0)
(45,0)
(55,0)
(65,0)
(75,0)
(85,0)
(95,0)
(105, 0.0)
};

\addplot[mark=none, red ,dashed, thick] coordinates {(14.512932376225978,0) (14.512932376225978,194.0)};
\addplot [red,, nodes near coords={\small Mean: 14.51\%},every node near coord/.style={anchor=180}] coordinates {( -100, 58.199999999999996)};
\addplot[mark=none, black ,dashed ,thick] coordinates {(13.923798785395924,0) (13.923798785395924,194.0)};
\addplot [black,, nodes near coords={\small Median: 13.92\%},every node near coord/.style={anchor=180}] coordinates {(-100, 53.35)};
\end{axis}
\end{tikzpicture}

%% file: figs/figure_hist_rhos.tex
\begin{figure}
\centering
\begin{subfigure}{.5\textwidth}
  \centering
    \captionsetup{width=.9\textwidth}
  \input{./figs/histograms/hist_rhos.tex}
  \label{fig:sub1}
    \subcaption{Reduction in density error by neural network enhancement of the algorithm when neural network was applied at a randomly chosen iteration. 
    Improvement was achieved in 60.0\% of cases.
    }
\end{subfigure}%
\begin{subfigure}{.5\textwidth}
  \centering
  \captionsetup{width=.9\textwidth}
  \input{./figs/histograms/hist_rhos_enhanced.tex}  
  \label{fig:sub2}
  \subcaption{Reduction in density error by neural network enhancement of the algorithm when neural network was applied via the previously described acceleration strategy.
  Improvement was achieved in 94.7\% of cases.
  }
\end{subfigure}
    \caption{Histograms of the reduction in density error by the neural network-accelerated algorithms, evaluated on the 488 relevant test cases. 
    In blue are the cases, where the enhanced and classical algorithm converged to the same local minimum, in orange the cases of different minima. 
    The dashed bars denote the cases with increase bigger than 100\%.
    }
\label{fig:figure_hist_rhos}
\end{figure}
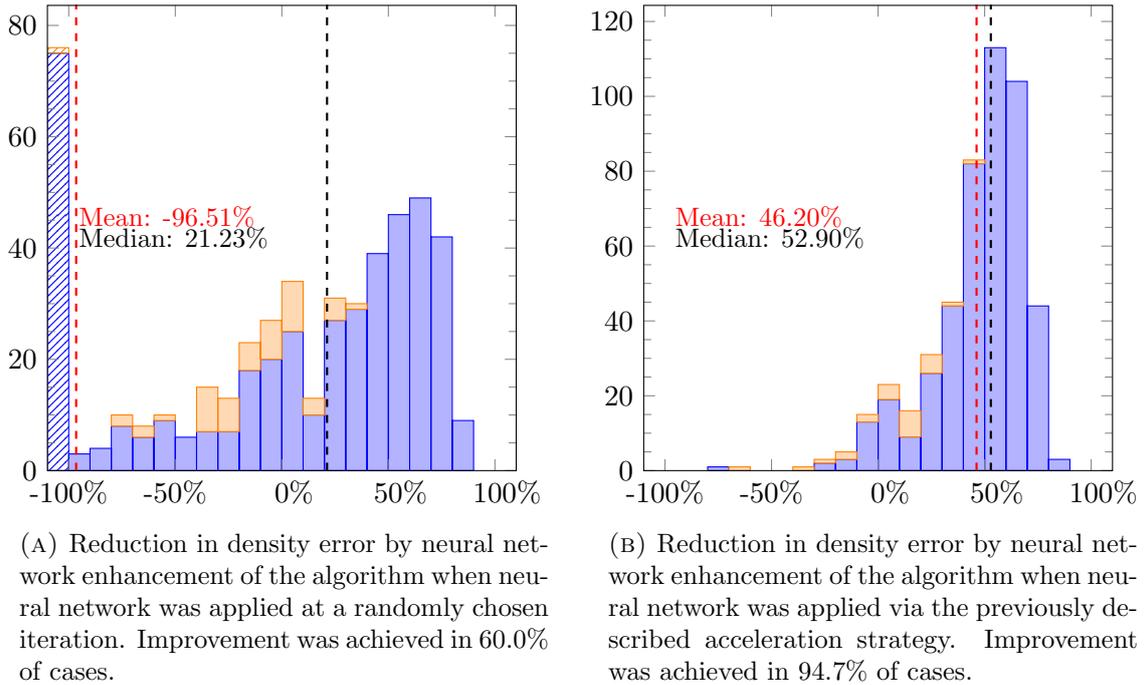

%% file: figs/histograms/hist_rhos.tex
\begin{tikzpicture}
\begin{axis}[
    width = \textwidth,
    height = \textwidth,
    ymin=0, ymax=83.60000000000001,
    xmin = -110, xmax = 110,
    minor y tick num = 3,
    xtick = {-100, -50, 0, 50, 100},
    xticklabels = {{-100\%, -50\%, \phantom{\%}0\%, \phantom{\%}50\%, \phantom{\%}100\%}},
    area style,
    ]
\addplot+[ybar stacked, bar width=10, color = blue!30, draw = blue] coordinates { 
(-105, 0)
(-95,3.0)
(-85,4.0)
(-75,8.0)
(-65,6.0)
(-55,9.0)
(-45,6.0)
(-35,7.0)
(-25,7.0)
(-15,18.0)
(-5,20.0)
(5,25.0)
(15,10.0)
(25,27.0)
(35,29.0)
(45,39.0)
(55,46.0)
(65,49.0)
(75,42.0)
(85,9.0)
(95,0.0)
(105, 0)
};
\addplot+[ybar stacked, bar width=10, color = blue!30, draw = blue, pattern = north east lines, pattern color = blue] coordinates { 
(-105, 75.0)
(-95,0)
(-85,0)
(-75,0)
(-65,0)
(-55,0)
(-45,0)
(-35,0)
(-25,0)
(-15,0)
(-5,0)
(5,0)
(15,0)
(25,0)
(35,0)
(45,0)
(55,0)
(65,0)
(75,0)
(85,0)
(95,0)
(105, 0.0)
};
\addplot+[ybar stacked, bar width=10, color = orange!30, draw = orange] coordinates { 
(-105, 0)
(-95,0.0)
(-85,0.0)
(-75,2.0)
(-65,2.0)
(-55,1.0)
(-45,0.0)
(-35,8.0)
(-25,6.0)
(-15,5.0)
(-5,7.0)
(5,9.0)
(15,3.0)
(25,4.0)
(35,1.0)
(45,0.0)
(55,0.0)
(65,0.0)
(75,0.0)
(85,0.0)
(95,0.0)
(105, 0)
};
\addplot+[ybar stacked, bar width=10, color = orange!30, draw = orange, pattern = north east lines, pattern color = orange] coordinates { 
(-105, 1.0)
(-95,0)
(-85,0)
(-75,0)
(-65,0)
(-55,0)
(-45,0)
(-35,0)
(-25,0)
(-15,0)
(-5,0)
(5,0)
(15,0)
(25,0)
(35,0)
(45,0)
(55,0)
(65,0)
(75,0)
(85,0)
(95,0)
(105, 0.0)
};

\addplot[mark=none, red ,dashed, thick] coordinates {(-96.51196273592873,0) (-96.51196273592873,152.0)};
\addplot [red,, nodes near coords={\small Mean: -96.51\%},every node near coord/.style={anchor=180}] coordinates {( -100, 45.6)};
\addplot[mark=none, black ,dashed ,thick] coordinates {(21.22814805440671,0) (21.22814805440671,152.0)};
\addplot [black,, nodes near coords={\small Median: 21.23\%},every node near coord/.style={anchor=180}] coordinates {(-100, 41.800000000000004)};
\end{axis}
\end{tikzpicture}

%% file: figs/histograms/hist_rhos_enhanced.tex
\begin{tikzpicture}
\begin{axis}[
    width = \textwidth,
    height = \textwidth,
    ymin=0, ymax=124.30000000000001,
    xmin = -110, xmax = 110,
    minor y tick num = 3,
    xtick = {-100, -50, 0, 50, 100},
    xticklabels = {{-100\%, -50\%, \phantom{\%}0\%, \phantom{\%}50\%, \phantom{\%}100\%}},
    area style,
    ]
\addplot+[ybar stacked, bar width=10, color = blue!30, draw = blue] coordinates { 
(-105, 0)
(-95,0.0)
(-85,0.0)
(-75,1.0)
(-65,0.0)
(-55,0.0)
(-45,0.0)
(-35,0.0)
(-25,2.0)
(-15,3.0)
(-5,13.0)
(5,19.0)
(15,9.0)
(25,26.0)
(35,44.0)
(45,82.0)
(55,113.0)
(65,104.0)
(75,44.0)
(85,3.0)
(95,0.0)
(105, 0)
};
\addplot+[ybar stacked, bar width=10, color = blue!30, draw = blue, pattern = north east lines, pattern color = blue] coordinates { 
(-105, 0.0)
(-95,0)
(-85,0)
(-75,0)
(-65,0)
(-55,0)
(-45,0)
(-35,0)
(-25,0)
(-15,0)
(-5,0)
(5,0)
(15,0)
(25,0)
(35,0)
(45,0)
(55,0)
(65,0)
(75,0)
(85,0)
(95,0)
(105, 0.0)
};
\addplot+[ybar stacked, bar width=10, color = orange!30, draw = orange] coordinates { 
(-105, 0)
(-95,0.0)
(-85,0.0)
(-75,0.0)
(-65,1.0)
(-55,0.0)
(-45,0.0)
(-35,1.0)
(-25,1.0)
(-15,2.0)
(-5,2.0)
(5,4.0)
(15,7.0)
(25,5.0)
(35,1.0)
(45,1.0)
(55,0.0)
(65,0.0)
(75,0.0)
(85,0.0)
(95,0.0)
(105, 0)
};
\addplot+[ybar stacked, bar width=10, color = orange!30, draw = orange, pattern = north east lines, pattern color = orange] coordinates { 
(-105, 0.0)
(-95,0)
(-85,0)
(-75,0)
(-65,0)
(-55,0)
(-45,0)
(-35,0)
(-25,0)
(-15,0)
(-5,0)
(5,0)
(15,0)
(25,0)
(35,0)
(45,0)
(55,0)
(65,0)
(75,0)
(85,0)
(95,0)
(105, 0.0)
};

\addplot[mark=none, red ,dashed, thick] coordinates {(46.19669297977423,0) (46.19669297977423,226.0)};
\addplot [red,, nodes near coords={\small Mean: 46.20\%},every node near coord/.style={anchor=180}] coordinates {( -100, 67.8)};
\addplot[mark=none, black ,dashed ,thick] coordinates {(52.896765979964115,0) (52.896765979964115,226.0)};
\addplot [black,, nodes near coords={\small Median: 52.90\%},every node near coord/.style={anchor=180}] coordinates {(-100, 62.150000000000006)};
\end{axis}
\end{tikzpicture}

%% file: figs/nn-best-rho-1.tex
\begin{figure}[H]
    \centering
    \begin{tabular}{ccc}
        \includegraphics[width = 0.3\textwidth, trim={3.8cm 0.2cm 3.8cm 0.8cm},clip]{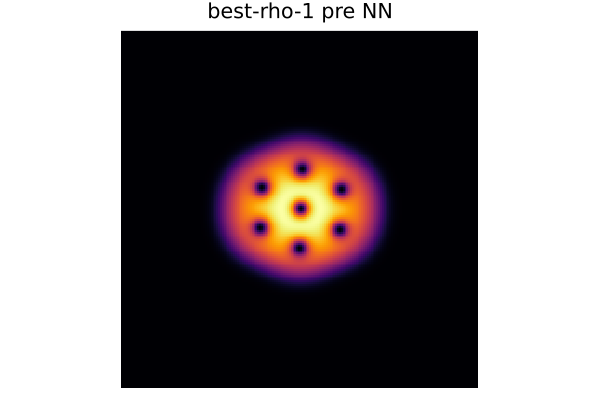} & 
        \includegraphics[width = 0.3\textwidth, trim={3.8cm 0.2cm 3.8cm 0.8cm},clip]{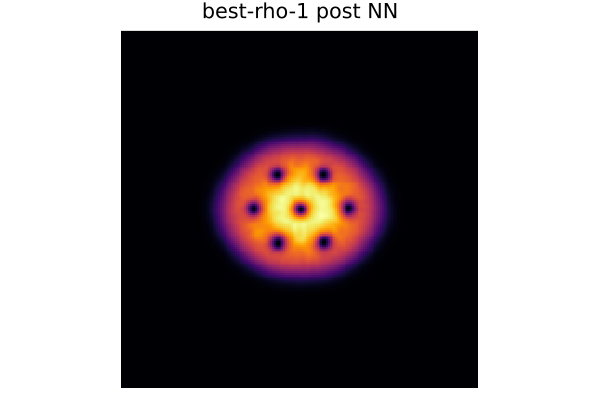} & 
        \includegraphics[width = 0.3\textwidth, trim={3.8cm 0.2cm 3.8cm 0.8cm},clip]{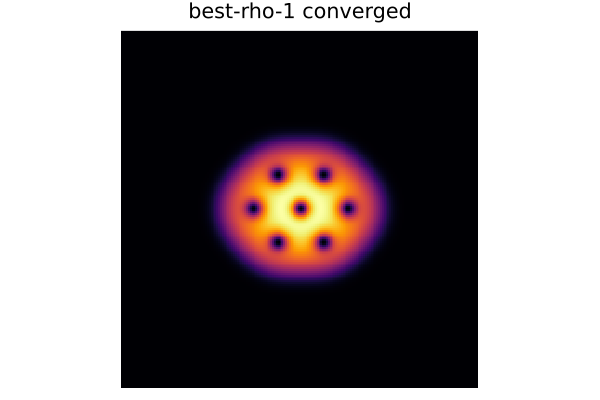} 
         \\
        NN in & NN out & post-NN result 
    \end{tabular}
    \caption{Neural network enhanced convergence of energy-adaptive RCG for parameters $v \approx (1.212, 1)$, $\omega \approx 1.396$,
    $\kappa = 538$. Enhancement resulted in 1046 total steps,
    while the classical algorithm took 1313 steps. Thus, the amount of steps was reduced by 20.34\%. 
    The density improved by 84.8\%. }
    \label{fig:nn-best-rho-1}
\end{figure}

%% file: figs/nn-best-rho-2.tex
\begin{figure}[H]
    \centering
    \begin{tabular}{ccc}
        \includegraphics[width = 0.3\textwidth, trim={3.8cm 0.2cm 3.8cm 0.8cm},clip]{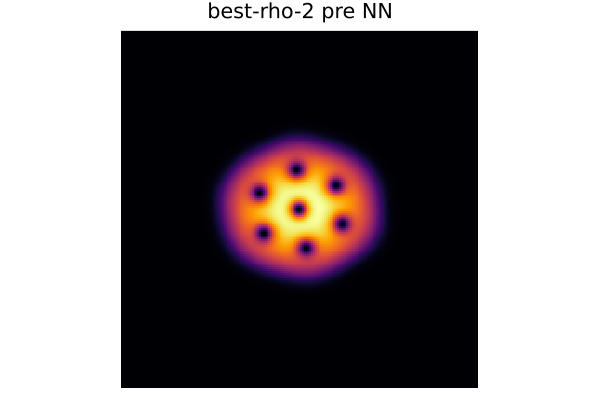} & 
        \includegraphics[width = 0.3\textwidth, trim={3.8cm 0.2cm 3.8cm 0.8cm},clip]{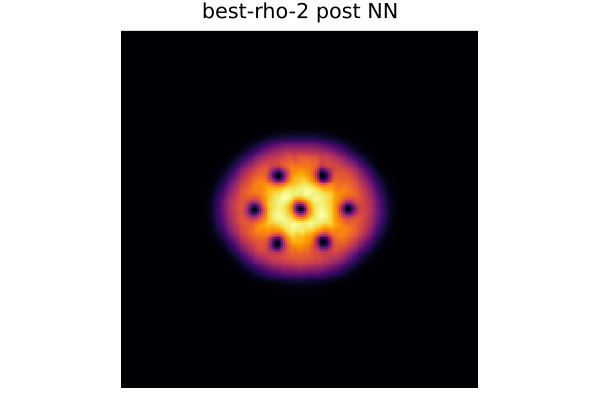} & 
        \includegraphics[width = 0.3\textwidth, trim={3.8cm 0.2cm 3.8cm 0.8cm},clip]{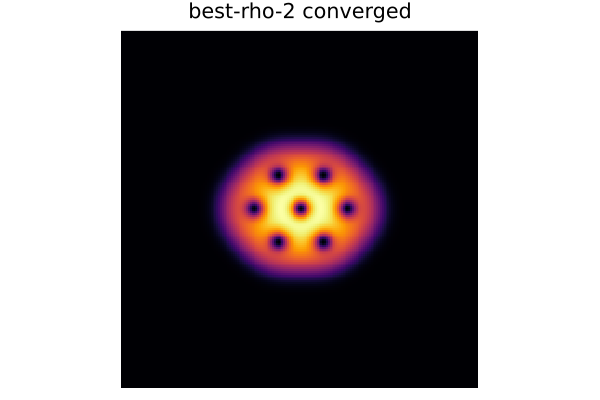} 
         \\
        NN in & NN out & post-NN result 
    \end{tabular}
    \caption{Neural network enhanced convergence of energy-adaptive RCG for parameters $v \approx (1.207, 1)$, $\omega \approx 1.431$,
    $\kappa = 477$. Enhancement resulted in 790 total steps,
    while the classical algorithm took 1047 steps. Thus, the amount of steps was reduced by 24.55\%. 
    The density improved by 82.9\%. }
    \label{fig:nn-best-rho-2}
\end{figure}

%% file: figs/nn-best-rho-3.tex
\begin{figure}[H]
    \centering
    \begin{tabular}{ccc}
        \includegraphics[width = 0.3\textwidth, trim={3.8cm 0.2cm 3.8cm 0.8cm},clip]{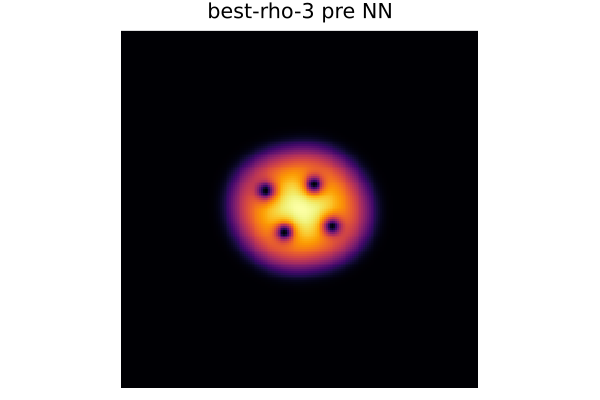} & 
        \includegraphics[width = 0.3\textwidth, trim={3.8cm 0.2cm 3.8cm 0.8cm},clip]{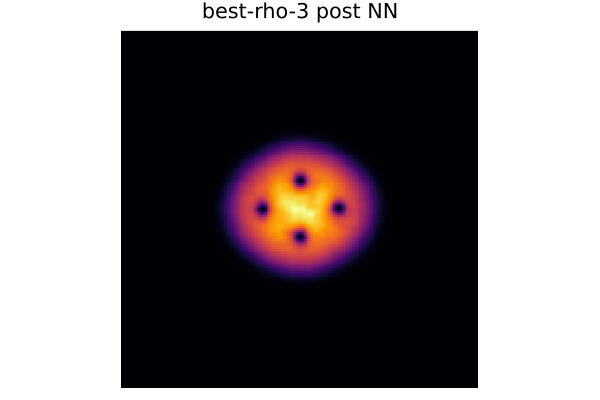} & 
        \includegraphics[width = 0.3\textwidth, trim={3.8cm 0.2cm 3.8cm 0.8cm},clip]{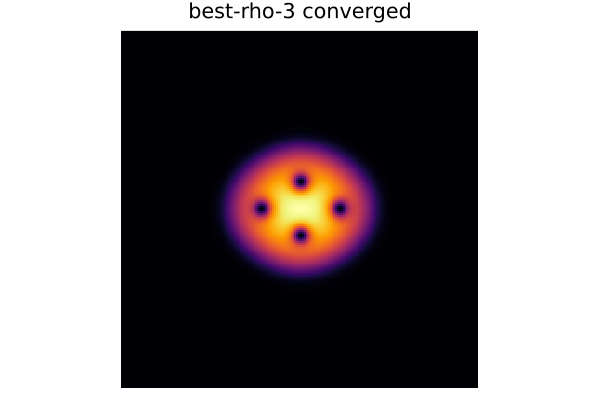} 
         \\
        NN in & NN out & post-NN result 
    \end{tabular}
    \caption{Neural network enhanced convergence of energy-adaptive RCG for parameters $v \approx (1.187, 1)$, $\omega \approx 1.115$,
    $\kappa = 442$. Enhancement resulted in 605 total steps,
    while the classical algorithm took 737 steps. Thus, the amount of steps was reduced by 17.91\%. 
    The density improved by 81.7\%. }
    \label{fig:nn-best-rho-3}
\end{figure}

%% file: figs/nn-worst-rho-1.tex
\begin{figure}[H]
    \centering
    \begin{tabular}{ccc}
        \includegraphics[width = 0.3\textwidth, trim={3.8cm 0.2cm 3.8cm 0.8cm},clip]{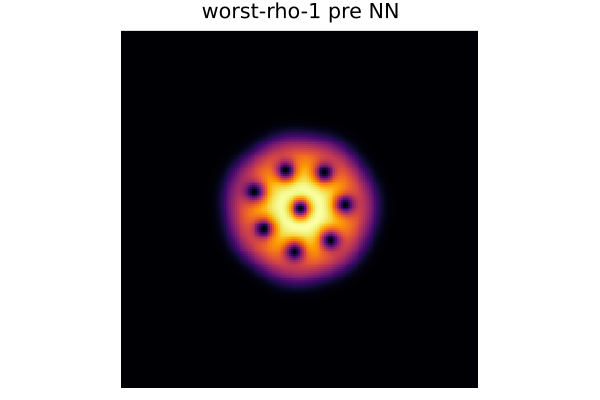} & 
        \includegraphics[width = 0.3\textwidth, trim={3.8cm 0.2cm 3.8cm 0.8cm},clip]{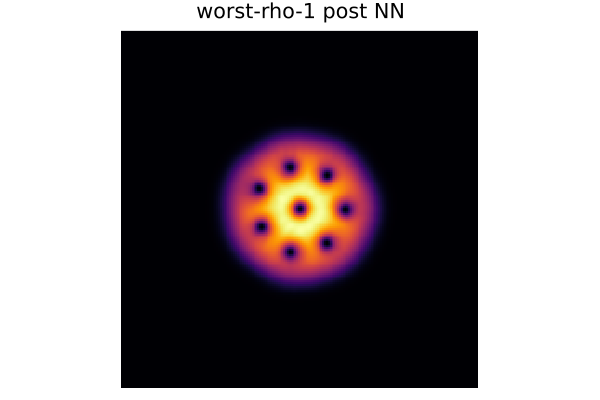} & 
        \includegraphics[width = 0.3\textwidth, trim={3.8cm 0.2cm 3.8cm 0.8cm},clip]{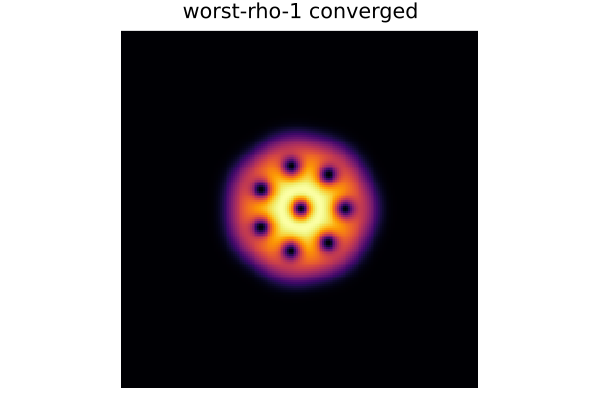} 
         \\
        NN in & NN out & post-NN result 
    \end{tabular}
    \caption{Neural network enhanced convergence of energy-adaptive RCG for parameters $v \approx (1.016, 1)$, $\omega \approx 1.485$,
    $\kappa = 357$. Enhancement resulted in 474 total steps,
    while the classical algorithm took 542 steps. Thus, the amount of steps was reduced by 12.55\%. 
    The density deteriorated by 79.1\%. }
    \label{fig:nn-worst-rho-1}
\end{figure}

%% file: figs/nn-worst-rho-2.tex
\begin{figure}[H]
    \centering
    \begin{tabular}{cccc}
        \includegraphics[width = 0.22\textwidth, trim={3.8cm 0.2cm 3.8cm 0.8cm},clip]{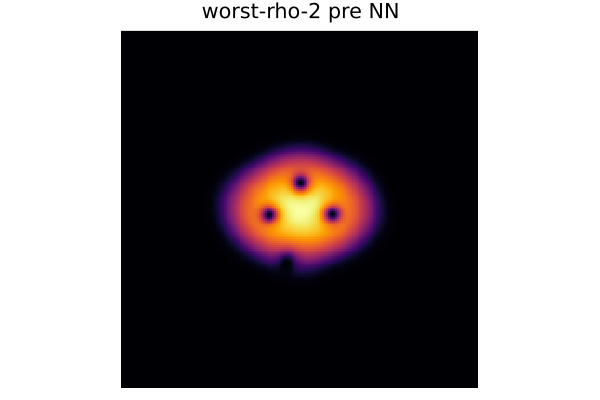} & 
        \includegraphics[width = 0.22\textwidth, trim={3.8cm 0.2cm 3.8cm 0.8cm},clip]{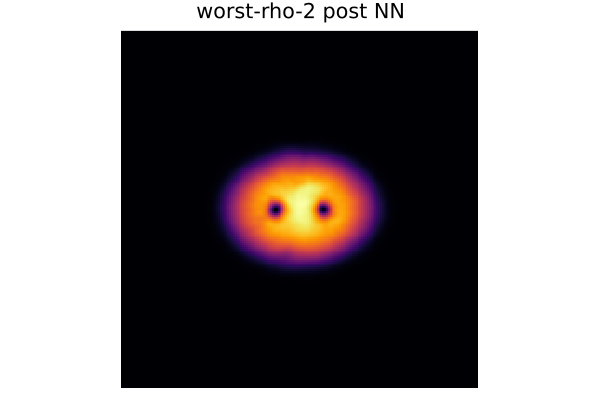} & 
        \includegraphics[width = 0.22\textwidth, trim={3.8cm 0.2cm 3.8cm 0.8cm},clip]{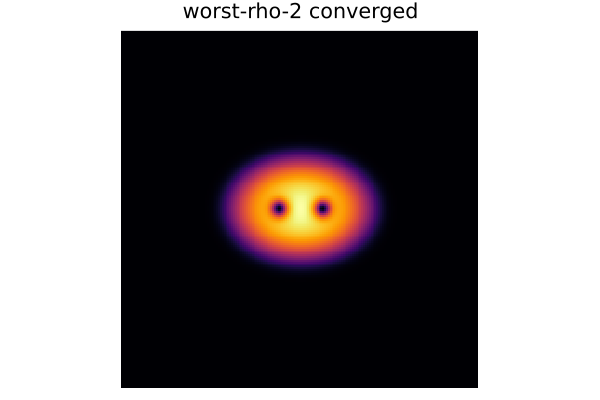} 
        & \includegraphics[width = 0.22\textwidth, trim={3.8cm 0.2cm 3.8cm 0.8cm},clip]{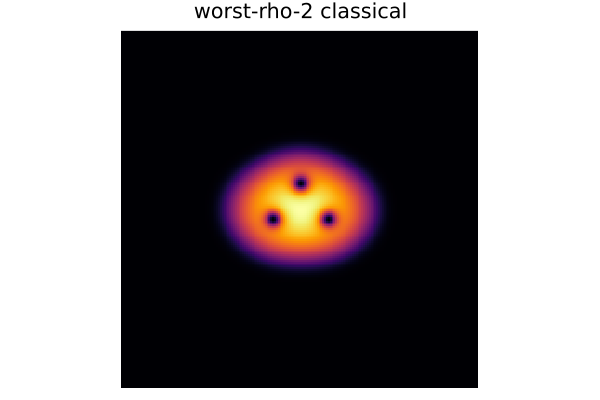} \\
        NN in & NN out & post-NN result & EARCG result
    \end{tabular}
    \caption{Neural network enhanced convergence of energy-adaptive RCG for parameters $v \approx (1.509, 1)$, $\omega \approx 1.111$,
    $\kappa = 468$. Enhancement resulted in 216 total steps,
    while the classical algorithm took 562 steps. Thus, the amount of steps was reduced by 61.57\%. 
    The density deteriorated by 67.7\%. The energy of the enhanced iteration was lower than the classical energy.}
    \label{fig:nn-worst-rho-2}
\end{figure}

%% file: figs/nn-worst-rho-3.tex
\begin{figure}[H]
    \centering
    \begin{tabular}{cccc}
        \includegraphics[width = 0.22\textwidth, trim={3.8cm 0.2cm 3.8cm 0.8cm},clip]{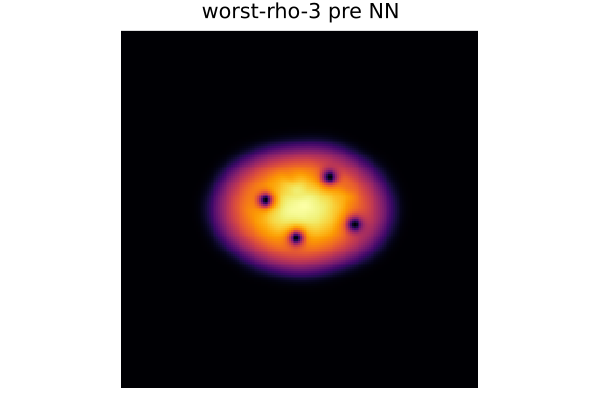} & 
        \includegraphics[width = 0.22\textwidth, trim={3.8cm 0.2cm 3.8cm 0.8cm},clip]{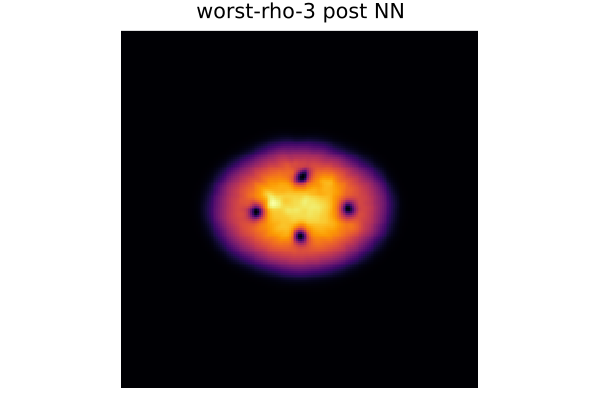} & 
        \includegraphics[width = 0.22\textwidth, trim={3.8cm 0.2cm 3.8cm 0.8cm},clip]{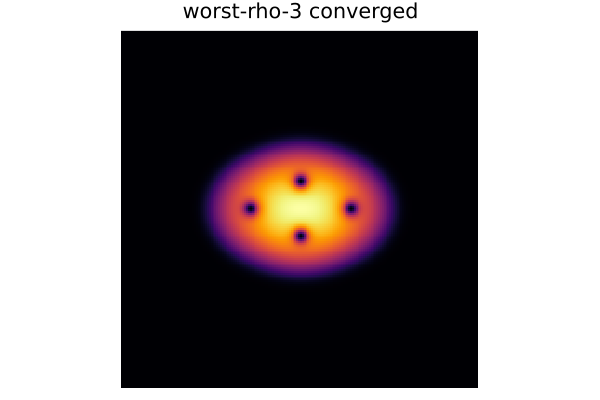} 
        & \includegraphics[width = 0.22\textwidth, trim={3.8cm 0.2cm 3.8cm 0.8cm},clip]{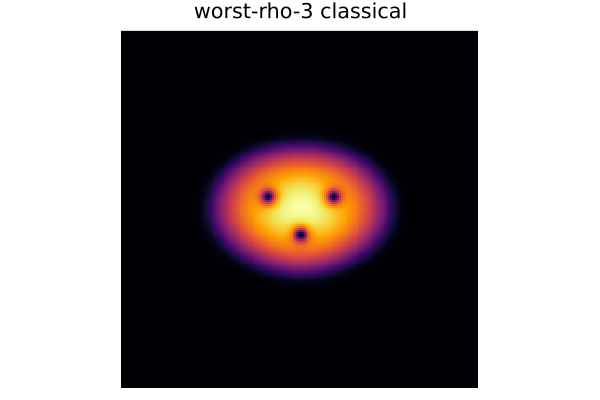} \\
        NN in & NN out & post-NN result & EARCG result
    \end{tabular}
    \caption{Neural network enhanced convergence of energy-adaptive RCG for parameters $v \approx (1.711, 1)$, $\omega \approx 0.896$,
    $\kappa = 938$. Enhancement resulted in 911 total steps,
    while the classical algorithm took 1005 steps. Thus, the amount of steps was reduced by 9.35\%. 
    The density deteriorated by 34.6\%. The energy of the enhanced iteration was higher than the classical energy.}
    \label{fig:nn-worst-rho-3}
\end{figure}

%% file: figs/nn-best-it-1.tex
\begin{figure}[H]
    \centering
    \begin{tabular}{ccc}
        \includegraphics[width = 0.3\textwidth, trim={3.8cm 0.2cm 3.8cm 0.8cm},clip]{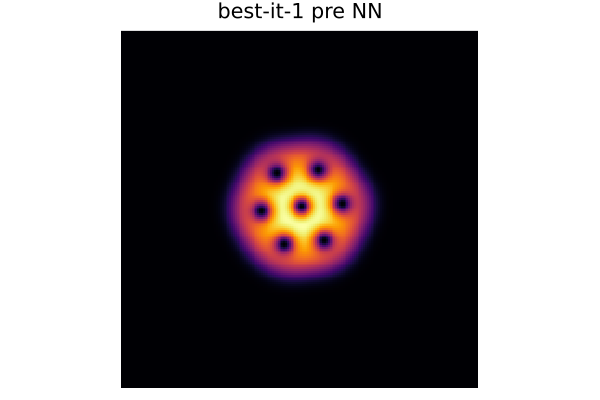} & 
        \includegraphics[width = 0.3\textwidth, trim={3.8cm 0.2cm 3.8cm 0.8cm},clip]{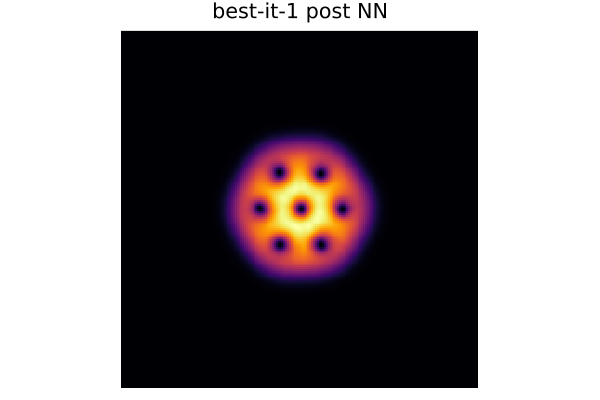} & 
        \includegraphics[width = 0.3\textwidth, trim={3.8cm 0.2cm 3.8cm 0.8cm},clip]{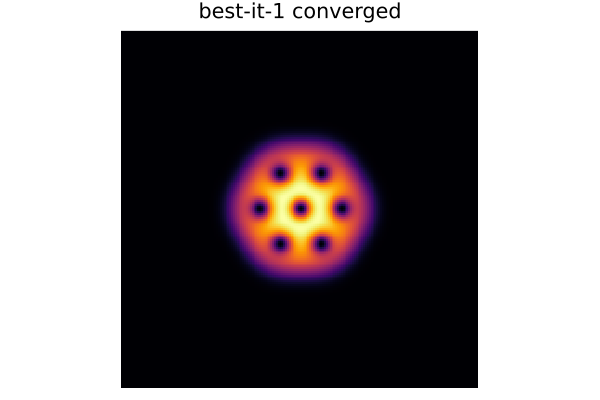} 
         \\
        NN in & NN out & post-NN result 
    \end{tabular}
    \caption{Neural network enhanced convergence of energy-adaptive RCG for parameters $v \approx (1.011, 1)$, $\omega \approx 1.524$,
    $\kappa = 271$. Enhancement resulted in 267 total steps,
    while the classical algorithm took 25285 steps. Thus, the amount of steps was reduced by 98.94\%. 
    The density improved by 62.8\%. }
    \label{fig:nn-best-it-1}
\end{figure}

%% file: figs/nn-best-it-2.tex
\begin{figure}[H]
    \centering
    \begin{tabular}{ccc}
        \includegraphics[width = 0.3\textwidth, trim={3.8cm 0.2cm 3.8cm 0.8cm},clip]{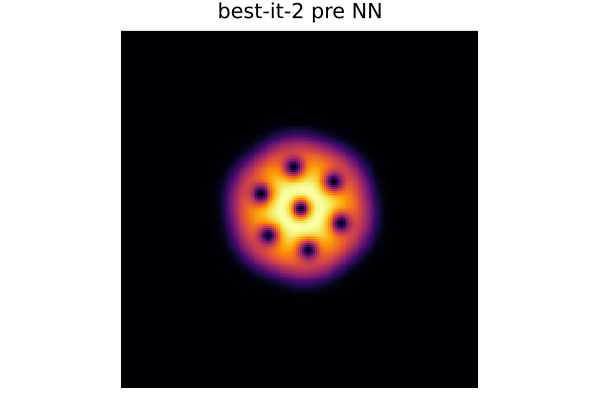} & 
        \includegraphics[width = 0.3\textwidth, trim={3.8cm 0.2cm 3.8cm 0.8cm},clip]{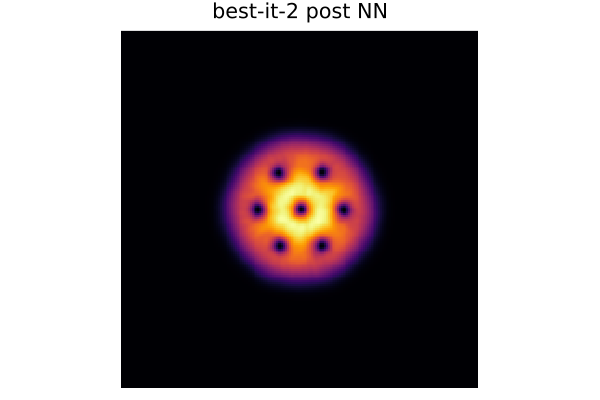} & 
        \includegraphics[width = 0.3\textwidth, trim={3.8cm 0.2cm 3.8cm 0.8cm},clip]{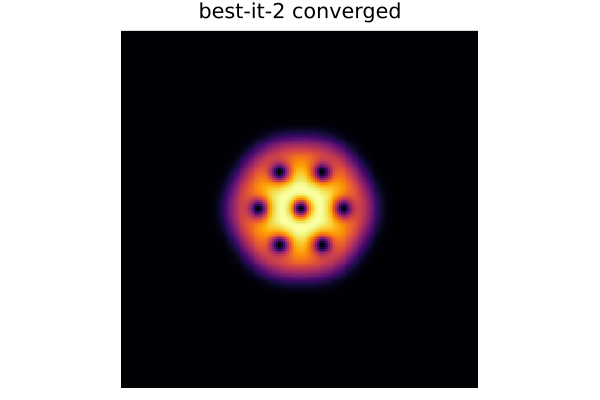} 
         \\
        NN in & NN out & post-NN result 
    \end{tabular}
    \caption{Neural network enhanced convergence of energy-adaptive RCG for parameters $v \approx (1.020, 1)$, $\omega \approx 1.423$,
    $\kappa = 385$. Enhancement resulted in 346 total steps,
    while the classical algorithm took 26206 steps. Thus, the amount of steps was reduced by 98.68\%. 
    The density improved by 78.8\%. }
    \label{fig:nn-best-it-2}
\end{figure}

%% file: figs/nn-best-it-3.tex
\begin{figure}[H]
    \centering
    \begin{tabular}{cccc}
        \includegraphics[width = 0.22\textwidth, trim={3.8cm 0.2cm 3.8cm 0.8cm},clip]{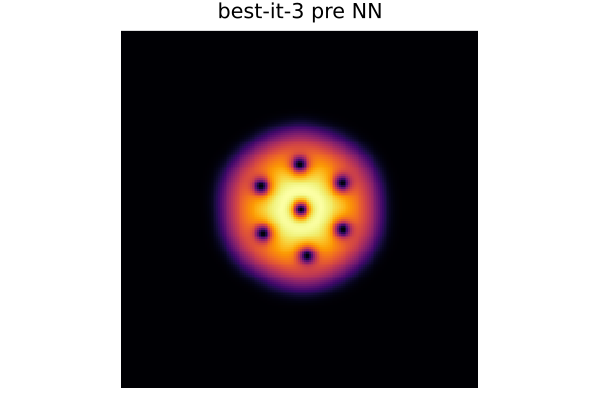} & 
        \includegraphics[width = 0.22\textwidth, trim={3.8cm 0.2cm 3.8cm 0.8cm},clip]{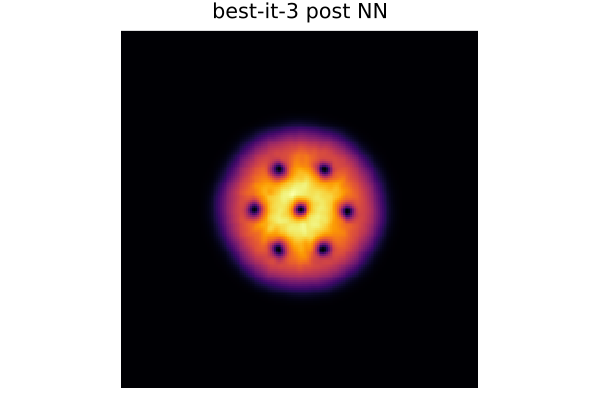} & 
        \includegraphics[width = 0.22\textwidth, trim={3.8cm 0.2cm 3.8cm 0.8cm},clip]{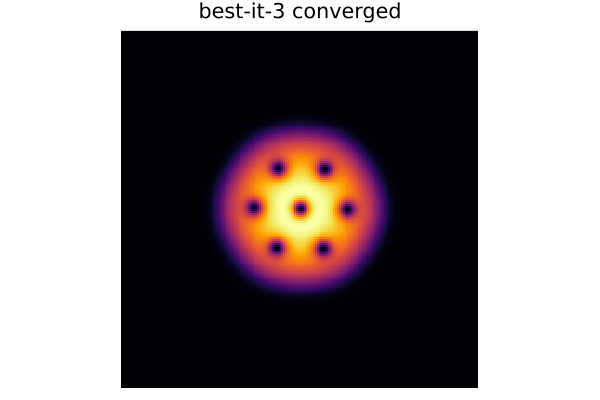} 
        & \includegraphics[width = 0.22\textwidth, trim={3.8cm 0.2cm 3.8cm 0.8cm},clip]{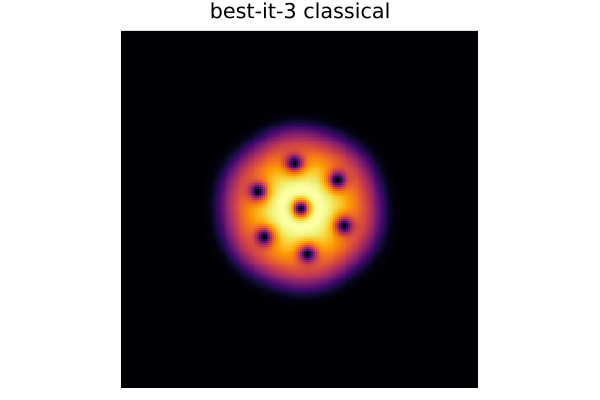} \\
        NN in & NN out & post-NN result & EARCG result
    \end{tabular}
    \caption{Neural network enhanced convergence of energy-adaptive RCG for parameters $v \approx (1.025, 1)$, $\omega \approx 1.207$,
    $\kappa = 725$. Enhancement resulted in 485 total steps,
    while the classical algorithm did not converge within 30000 steps. Thus, the amount of steps was reduced by 98.38\%. 
    Since the classical algorithm did not converge, no estimation of the density improvement can be given. The energy of the enhanced iteration was lower than the classical energy.}
    \label{fig:nn-best-it-3}
\end{figure}

%% file: figs/nn-worst-it-1.tex
\begin{figure}[H]
    \centering
    \begin{tabular}{ccc}
        \includegraphics[width = 0.3\textwidth, trim={3.8cm 0.2cm 3.8cm 0.8cm},clip]{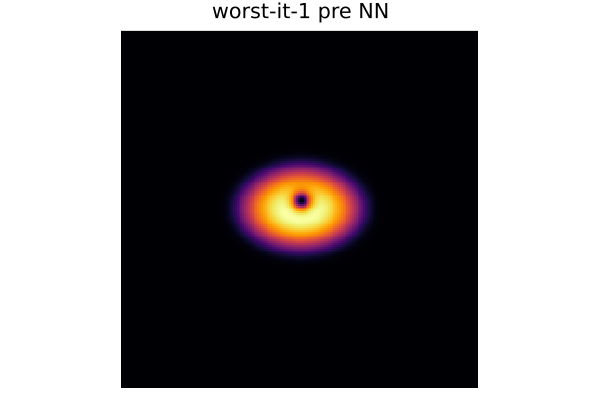} & 
        \includegraphics[width = 0.3\textwidth, trim={3.8cm 0.2cm 3.8cm 0.8cm},clip]{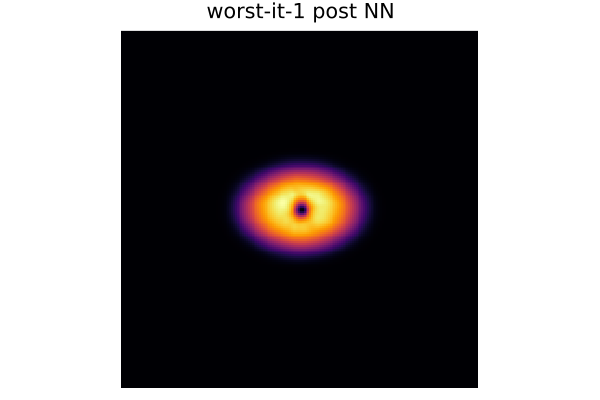} & 
        \includegraphics[width = 0.3\textwidth, trim={3.8cm 0.2cm 3.8cm 0.8cm},clip]{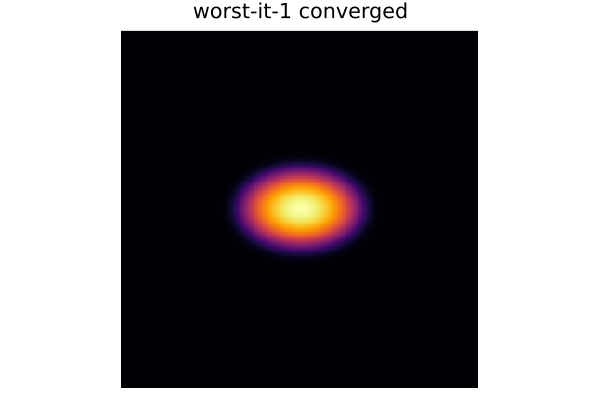} 
         \\
        NN in & NN out & post-NN result 
    \end{tabular}
    \caption{Neural network enhanced convergence of energy-adaptive RCG for parameters $v \approx (1.970, 1)$, $\omega \approx 0.829$,
    $\kappa = 248$. Enhancement resulted in 797 total steps,
    while the classical algorithm took 577 steps. Thus, the amount of steps was increased by 38.13\%. 
    The density deteriorated by 2.7\%. }
    \label{fig:nn-worst-it-1}
\end{figure}

%% file: figs/nn-worst-it-2.tex
\begin{figure}[H]
    \centering
    \begin{tabular}{ccc}
        \includegraphics[width = 0.3\textwidth, trim={3.8cm 0.2cm 3.8cm 0.8cm},clip]{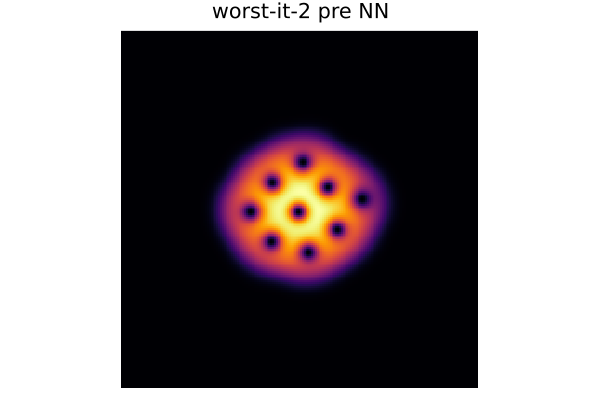} & 
        \includegraphics[width = 0.3\textwidth, trim={3.8cm 0.2cm 3.8cm 0.8cm},clip]{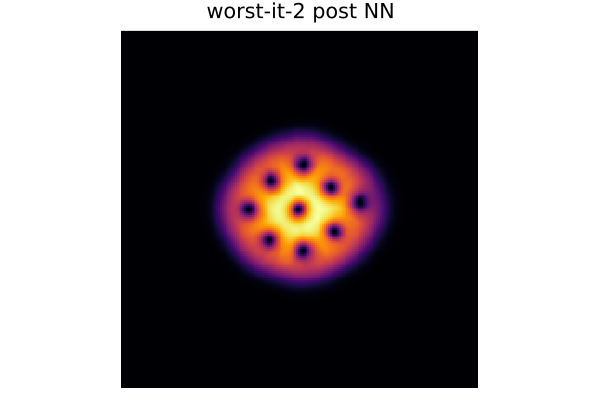} & 
        \includegraphics[width = 0.3\textwidth, trim={3.8cm 0.2cm 3.8cm 0.8cm},clip]{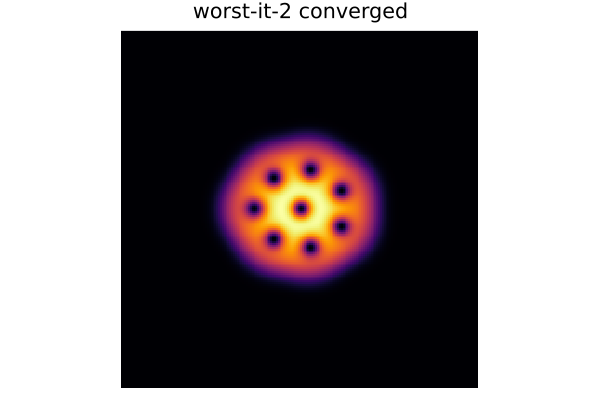} 
         \\
        NN in & NN out & post-NN result 
    \end{tabular}
    \caption{Neural network enhanced convergence of energy-adaptive RCG for parameters $v \approx (1.092, 1)$, $\omega \approx 1.563$,
    $\kappa = 393$. Enhancement resulted in 4525 total steps,
    while the classical algorithm took 3530 steps. Thus, the amount of steps was increased by 28.19\%. 
    The density improved by 0.8\%. }
    \label{fig:nn-worst-it-2}
\end{figure}

%% file: figs/nn-worst-it-3.tex
\begin{figure}[H]
    \centering
    \begin{tabular}{cccc}
        \includegraphics[width = 0.22\textwidth, trim={3.8cm 0.2cm 3.8cm 0.8cm},clip]{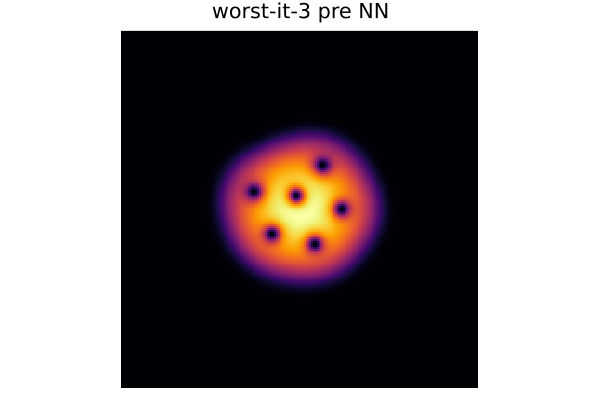} & 
        \includegraphics[width = 0.22\textwidth, trim={3.8cm 0.2cm 3.8cm 0.8cm},clip]{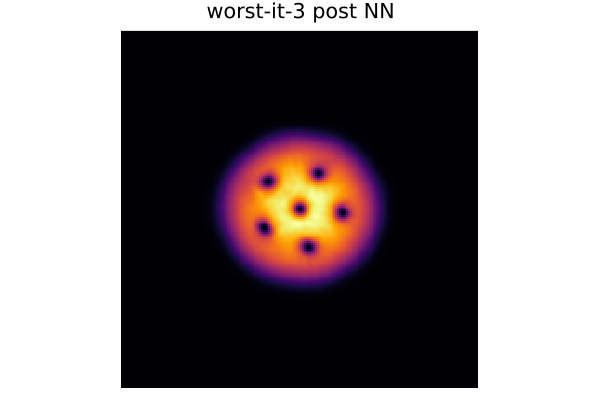} & 
        \includegraphics[width = 0.22\textwidth, trim={3.8cm 0.2cm 3.8cm 0.8cm},clip]{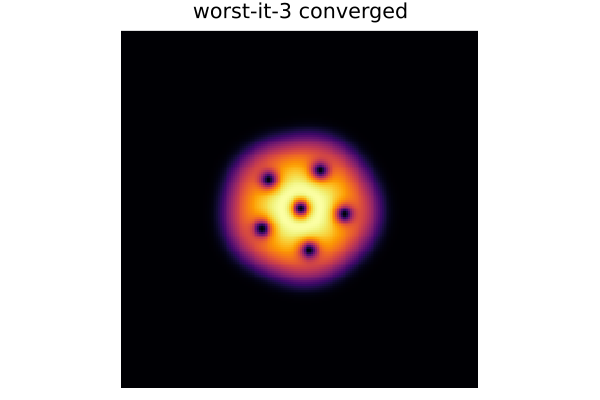} 
        & \includegraphics[width = 0.22\textwidth, trim={3.8cm 0.2cm 3.8cm 0.8cm},clip]{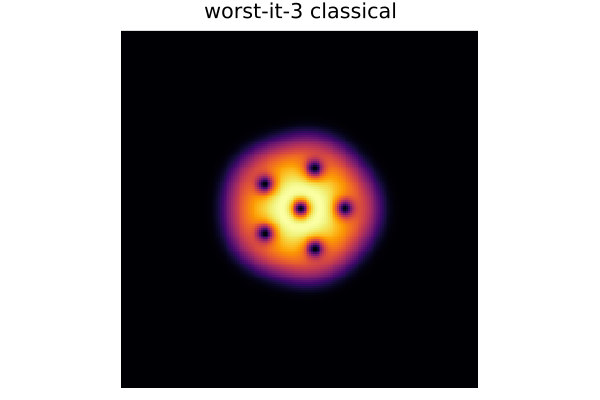} \\
        NN in & NN out & post-NN result & EARCG result
    \end{tabular}
    \caption{Neural network enhanced convergence of energy-adaptive RCG for parameters $v \approx (1.082, 1)$, $\omega \approx 1.219$,
    $\kappa = 599$. Enhancement resulted in no convergence within 30000 total steps,
    while the classical algorithm took 24068 steps. Thus, the amount of steps was increased by 24.65\%. 
    The density improved by 30.1\%. The energy of the enhanced iteration was lower than the classical energy.}
    \label{fig:nn-worst-it-3}
\end{figure}

%% file: main.bbl
\begin{thebibliography}{10}

\bibitem{adler2017solving}
J.~Adler and O.~{\"O}ktem.
\newblock Solving ill-posed inverse problems using iterative deep neural
  networks.
\newblock {\em Inverse Problems}, 33(12):124007, 2017.

\bibitem{AltHP21}
R.~Altmann, P.~Henning, and D.~Peterseim.
\newblock The {$J$}-method for the {G}ross--{P}itaevskii eigenvalue problem.
\newblock {\em Numer. Math.}, 148:575--610, 2021.

\bibitem{altmann2025riemannianoptimisationmethodsground}
R.~Altmann, M.~Hermann, D.~Peterseim, and T.~Stykel.
\newblock Riemannian optimisation methods for ground states of multicomponent
  bose-einstein condensates, 2025.

\bibitem{AltPS22}
R.~Altmann, D.~Peterseim, and T.~Stykel.
\newblock Energy-adaptive {R}iemannian optimization on the {S}tiefel manifold.
\newblock {\em ESAIM Math. Model. Numer. Anal.}, 56(5):1629--1653, 2022.

\bibitem{AltPS24}
R.~Altmann, D.~Peterseim, and T.~Stykel.
\newblock Riemannian {N}ewton methods for energy minimization problems of
  {K}ohn--{S}ham type.
\newblock {\em J. Sci. Comput.}, 101, 2024.

\bibitem{AntLT17}
X.~Antoine, A.~Levitt, and Q.~Tang.
\newblock Efficient spectral computation of the stationary states of rotating
  {B}ose–{E}instein condensates by preconditioned nonlinear conjugate
  gradient methods.
\newblock {\em J. Comput. Phys.}, 343:92--109, 2017.

\bibitem{PhysRevResearch.7.013332}
X.-D. Bai, T.~Xu, J.~Li, Y.-K. Liu, Y.~Zhao, and J.~Zhao.
\newblock Rapid discovering ground states in lee-huang-yang spin-orbit coupled
  bose-einstein condensates via a coupled-tgnn surrogate model.
\newblock {\em Phys. Rev. Res.}, 7:013332, Mar 2025.

\bibitem{BaoCai2013}
W.~Bao and Y.~Cai.
\newblock Mathematical theory and numerical methods for bose-einstein
  condensation.
\newblock {\em Kinetic and Related Models}, 6(1):1--135, 2013.

\bibitem{BAO2025113486}
W.~Bao, Z.~Chang, and X.~Zhao.
\newblock Computing ground states of {B}ose-{E}instein condensation by
  normalized deep neural network.
\newblock {\em Journal of Computational Physics}, 520:113486, 2025.

\bibitem{BaoCL06}
W.~Bao, I-L. Chern, and F.~Y. Lim.
\newblock Efficient and spectrally accurate numerical methods for computing
  ground and first excited states in {B}ose--{E}instein condensates.
\newblock {\em J. Comput. Phys.}, 219(2):836--854, 2006.

\bibitem{BaoD04}
W.~Bao and Q.~Du.
\newblock Computing the ground state solution of {B}ose--{E}instein condensates
  by a~normalized gradient flow.
\newblock {\em SIAM J. Sci. Comput.}, 25(5):1674--1697, 2004.

\bibitem{BaoT03}
W.~Bao and W.~Tang.
\newblock Ground-state solution of {B}ose--{E}instein condensate by directly
  minimizing the energy functional.
\newblock {\em J. Comput. Phys.}, 187(1):230--254, 2003.

\bibitem{CanCM10}
E.~Canc{\`e}s, R.~Chakir, and Y.~Maday.
\newblock Numerical analysis of nonlinear eigenvalue problems.
\newblock {\em J. Sci. Comput.}, 45(1-3):90--117, 2010.

\bibitem{CanKL21}
E.~Canc\`{e}s, G.~Kemlin, and A.~Levitt.
\newblock Convergence analysis of direct minimization and self-consistent
  iterations.
\newblock {\em SIAM J. Matrix Anal. Appl.}, 42(1):243--274, 2021.

\bibitem{ChenLLZ24}
Z.~Chen, J.~Lu, Y.~Lu, and X.~Zhang.
\newblock On the convergence of {S}obolev gradient flow for the
  {G}ross--{P}itaevskii eigenvalue problem.
\newblock {\em SIAM J. Numer. Anal.}, 62(2):667--691, 2024.

\bibitem{DanK10}
I.~Danaila and P.~Kazemi.
\newblock A new {S}obolev gradient method for direct minimization of the
  {G}ross–{P}itaevskii energy with rotation.
\newblock {\em SIAM J. Sci. Comput.}, 32(5):2447--2467, 2010.

\bibitem{DanP17}
I.~Danaila and B.~Protas.
\newblock Computation of ground states of the {G}ross--{P}itaevskii functional
  via {R}iemannian optimization.
\newblock {\em SIAM J. Sci. Comput.}, 39(6):B1102--B1129, 2017.

\bibitem{DioC07}
C.~M. Dion and E.~Canc{\`e}s.
\newblock Ground state of the time-independent {G}ross--{P}itaevskii equation.
\newblock {\em Comput. Phys. Commun.}, 177(10):787--798, 2007.

\bibitem{DuL22}
C.-E. Du and C.-S. Liu.
\newblock Newton--{N}oda iteration for computing the ground states of nonlinear
  {S}chr{\"o}dinger equations.
\newblock {\em SIAM J. Sci. Comput.}, 44(4):A2370--A2385, 2022.

\bibitem{GarP01}
J.~J. Garc{\'{\i}}a-Ripoll and V.~M. P{\'e}rez-Garc{\'{\i}}a.
\newblock Optimizing {S}chr\"odinger functionals using {S}obolev gradients:
  applications to quantum mechanics and nonlinear optics.
\newblock {\em SIAM J. Sci. Comput.}, 23(4):1316--1334, 2001.

\bibitem{HLP24}
M.~Hauck, Y.~Liang, and D.~Peterseim.
\newblock Positivity preserving finite element method for the gross-pitaevskii
  ground state: discrete uniqueness and global convergence, 2024.

\bibitem{HenJ23}
P.~Henning and E.~Jarlebring.
\newblock The gross–pitaevskii equation and eigenvector nonlinearities:
  Numerical methods and algorithms.
\newblock {\em SIAM Review}, 67(2):256--317, 2025.

\bibitem{HenP20}
P.~Henning and D.~Peterseim.
\newblock Sobolev gradient flow for the {G}ross--{P}itaevskii eigenvalue
  problem: global convergence and computational efficiency.
\newblock {\em SIAM J. Numer. Anal.}, 58(3):1744--1772, 2020.

\bibitem{DFTKjcon}
M.~F. Herbst, A.~Levitt, and E.~Canc{\`e}s.
\newblock {DFTK}: {A} {J}ulian approach for simulating electrons in solids.
\newblock {\em Proc. JuliaCon Conf.}, 3:69, 2021.

\bibitem{JarKM14}
E.~Jarlebring, S.~Kvaal, and W.~Michiels.
\newblock An inverse iteration method for eigenvalue problems with eigenvector
  nonlinearities.
\newblock {\em SIAM J. Sci. Comput.}, 36(4):A1978--A2001, 2014.

\bibitem{KaE10}
P.~Kazemi and M.~Eckart.
\newblock Minimizing the {G}ross--{P}itaevskii energy functional with the
  {S}obolev gradient -- analytical and numerical results.
\newblock {\em Int. J. Comput. Methods}, 7(3):453--475, 2010.

\bibitem{KohS65}
W.~Kohn and L.~J. Sham.
\newblock Self-consistent equations including exchange and correlation effects.
\newblock {\em Phys. Rev.}, 140:A1133--A1138, 1965.

\bibitem{LieSY00}
E.~H. Lieb, R.~Seiringer, and J.~Yngvason.
\newblock Bosons in a trap: {A} rigorous derivation of the {G}ross-{P}itaevskii
  energy functional.
\newblock {\em Phys. Rev. A}, 61, 2000.

\bibitem{PhysRevLett.78.586}
C.~J. Myatt, E.~A. Burt, R.~W. Ghrist, E.~A. Cornell, and C.~E. Wieman.
\newblock Production of two overlapping {B}ose-{E}instein condensates by
  sympathetic cooling.
\newblock {\em Phys. Rev. Lett.}, 78:586--589, 1997.

\bibitem{PPS25}
D.~Peterseim, J.~P\"uschel, and T.~Stykel.
\newblock Energy-adaptive riemannian conjugate gradient method for density
  functional theory.
\newblock 2025.

\bibitem{unet2015}
O.~Ronneberger, P.~Fischer, and T.~Brox.
\newblock U-net: Convolutional networks for biomedical image segmentation,
  2015.

\bibitem{Slater1951}
J.~C. Slater.
\newblock A simplification of the hartree-fock method.
\newblock {\em Phys. Rev.}, 81:385--390, Feb 1951.

\bibitem{SteISMCK98}
J.~Stenger, S.~Inouye, D.~M. Stamper-Kurn, H.-J. Miesner, A.~P. Chikkatur, and
  W.~Ketterle.
\newblock Spin domains in ground-state {B}ose-{E}instein condensates.
\newblock {\em Nature}, 396:345–348, 1998.

\bibitem{DBLP:journals/corr/abs-2107-10254}
S.~Venkataraman and B.~Amos.
\newblock Neural fixed-point acceleration for convex optimization.
\newblock {\em CoRR}, abs/2107.10254, 2021.

\bibitem{WuWB17}
X.~Wu, Z.~Wen, and W.~Bao.
\newblock A regularized {N}ewton method for computing ground states of
  {B}ose–{E}instein condensates.
\newblock {\em J. Sci. Comput.}, 73:303--329, 2017.

\end{thebibliography}
